\documentclass[reqno,12pt,a4paper]{amsart}
\usepackage{mathrsfs}
\usepackage{amssymb}
\usepackage{amsfonts}
\usepackage{latexsym}
\usepackage{amsthm}
\usepackage{graphicx}

\def\lb{\label}
\newcommand{\er}[1]{\textrm{(\ref{#1})}}

\begin{document}
%%%%%%%%%% Some definitions %%%%%%%%%%
%%%%%%%% Equations, theorems %%%%%%%%%
\renewcommand{\theequation}{\arabic{section}.\arabic{equation}}
\theoremstyle{plain}
\newtheorem{theorem}{\bf Theorem}[section]
\newtheorem{lemma}[theorem]{\bf Lemma}
\newtheorem{corollary}[theorem]{\bf Corollary}
\newtheorem{proposition}[theorem]{\bf Proposition}
\newtheorem{definition}[theorem]{\bf Definition}
\newtheorem{remark}[theorem]{\bf Remark}
%%%%% Alphabet %%%%%
\def\a{\alpha}  \def\cA{{\mathcal A}}     \def\bA{{\bf A}}  \def\mA{{\mathscr A}}
\def\b{\beta}   \def\cB{{\mathcal B}}     \def\bB{{\bf B}}  \def\mB{{\mathscr B}}
\def\g{\gamma}  \def\cC{{\mathcal C}}     \def\bC{{\bf C}}  \def\mC{{\mathscr C}}
\def\G{\Gamma}  \def\cD{{\mathcal D}}     \def\bD{{\bf D}}  \def\mD{{\mathscr D}}
\def\d{\delta}  \def\cE{{\mathcal E}}     \def\bE{{\bf E}}  \def\mE{{\mathscr E}}
\def\D{\Delta}  \def\cF{{\mathcal F}}     \def\bF{{\bf F}}  \def\mF{{\mathscr F}}
\def\c{\chi}    \def\cG{{\mathcal G}}     \def\bG{{\bf G}}  \def\mG{{\mathscr G}}
\def\z{\zeta}   \def\cH{{\mathcal H}}     \def\bH{{\bf H}}  \def\mH{{\mathscr H}}
\def\e{\eta}    \def\cI{{\mathcal I}}     \def\bI{{\bf I}}  \def\mI{{\mathscr I}}
\def\p{\psi}    \def\cJ{{\mathcal J}}     \def\bJ{{\bf J}}  \def\mJ{{\mathscr J}}
\def\vT{\Theta} \def\cK{{\mathcal K}}     \def\bK{{\bf K}}  \def\mK{{\mathscr K}}
\def\k{\kappa}  \def\cL{{\mathcal L}}     \def\bL{{\bf L}}  \def\mL{{\mathscr L}}
\def\l{\lambda} \def\cM{{\mathcal M}}     \def\bM{{\bf M}}  \def\mM{{\mathscr M}}
\def\L{\Lambda} \def\cN{{\mathcal N}}     \def\bN{{\bf N}}  \def\mN{{\mathscr N}}
\def\m{\mu}     \def\cO{{\mathcal O}}     \def\bO{{\bf O}}  \def\mO{{\mathscr O}}
\def\n{\nu}     \def\cP{{\mathcal P}}     \def\bP{{\bf P}}  \def\mP{{\mathscr P}}
\def\r{\rho}    \def\cQ{{\mathcal Q}}     \def\bQ{{\bf Q}}  \def\mQ{{\mathscr Q}}
\def\s{\sigma}  \def\cR{{\mathcal R}}     \def\bR{{\bf R}}  \def\mR{{\mathscr R}}
                \def\cS{{\mathcal S}}     \def\bS{{\bf S}}  \def\mS{{\mathscr S}}
\def\t{\tau}    \def\cT{{\mathcal T}}     \def\bT{{\bf T}}  \def\mT{{\mathscr T}}
\def\f{\phi}    \def\cU{{\mathcal U}}     \def\bU{{\bf U}}  \def\mU{{\mathscr U}}
\def\F{\Phi}    \def\cV{{\mathcal V}}     \def\bV{{\bf V}}  \def\mV{{\mathscr V}}
\def\P{\Psi}    \def\cW{{\mathcal W}}     \def\bW{{\bf W}}  \def\mW{{\mathscr W}}
\def\o{\omega}  \def\cX{{\mathcal X}}     \def\bX{{\bf X}}  \def\mX{{\mathscr X}}
\def\x{\xi}     \def\cY{{\mathcal Y}}     \def\bY{{\bf Y}}  \def\mY{{\mathscr Y}}
\def\X{\Xi}     \def\cZ{{\mathcal Z}}     \def\bZ{{\bf Z}}  \def\mZ{{\mathscr Z}}
\def\O{\Omega}
\def\vs{\varsigma}

\def\mb{{\mathscr b}}
\def\mh{{\mathscr h}}
\def\me{{\mathscr e}}
\def\mk{{\mathscr k}}
\def\mz{{\mathscr z}}
\def\mx{{\mathscr x}}
\newcommand{\gA}{\mathfrak{A}}          \newcommand{\ga}{\mathfrak{a}}
\newcommand{\gB}{\mathfrak{B}}          \newcommand{\gb}{\mathfrak{b}}
\newcommand{\gC}{\mathfrak{C}}          \newcommand{\gc}{\mathfrak{c}}
\newcommand{\gD}{\mathfrak{D}}          \newcommand{\gd}{\mathfrak{d}}
\newcommand{\gE}{\mathfrak{E}}
\newcommand{\gF}{\mathfrak{F}}           \newcommand{\gf}{\mathfrak{f}}
\newcommand{\gG}{\mathfrak{G}}           %\newcommand{\gg}{\mathfrak{g}}
\newcommand{\gH}{\mathfrak{H}}           \newcommand{\gh}{\mathfrak{h}}
\newcommand{\gI}{\mathfrak{I}}           \newcommand{\gi}{\mathfrak{i}}
\newcommand{\gJ}{\mathfrak{J}}           \newcommand{\gj}{\mathfrak{j}}
\newcommand{\gK}{\mathfrak{K}}            \newcommand{\gk}{\mathfrak{k}}
\newcommand{\gL}{\mathfrak{L}}            \newcommand{\gl}{\mathfrak{l}}
\newcommand{\gM}{\mathfrak{M}}            \newcommand{\gm}{\mathfrak{m}}
\newcommand{\gN}{\mathfrak{N}}            \newcommand{\gn}{\mathfrak{n}}
\newcommand{\gO}{\mathfrak{O}}
\newcommand{\gP}{\mathfrak{P}}             \newcommand{\gp}{\mathfrak{p}}
\newcommand{\gQ}{\mathfrak{Q}}             \newcommand{\gq}{\mathfrak{q}}
\newcommand{\gR}{\mathfrak{R}}             \newcommand{\gr}{\mathfrak{r}}
\newcommand{\gS}{\mathfrak{S}}              \newcommand{\gs}{\mathfrak{s}}
\newcommand{\gT}{\mathfrak{T}}             \newcommand{\gt}{\mathfrak{t}}
\newcommand{\gU}{\mathfrak{U}}             \newcommand{\gu}{\mathfrak{u}}
\newcommand{\gV}{\mathfrak{V}}             \newcommand{\gv}{\mathfrak{v}}
\newcommand{\gW}{\mathfrak{W}}             \newcommand{\gw}{\mathfrak{w}}
\newcommand{\gX}{\mathfrak{X}}               \newcommand{\gx}{\mathfrak{x}}
\newcommand{\gY}{\mathfrak{Y}}              \newcommand{\gy}{\mathfrak{y}}
\newcommand{\gZ}{\mathfrak{Z}}             \newcommand{\gz}{\mathfrak{z}}

\def\ve{\varepsilon} \def\vt{\vartheta} \def\vp{\varphi}
\def\vk{\varkappa}
\def\vr{\varrho}

\def\be{{\bf e}} \def\bc{{\bf c}}
\def\bx{{\bf x}} \def\by{{\bf y}}
\def\bv{{\bf v}} \def\bu{{\bf u}}
 \def\bp{{\bf p}}
%************************
\def\mm{\mathrm m}
\def\mn{\mathrm n}

\def\A{{\mathbb A}}\def\B{{\mathbb B}}\def\C{{\mathbb C}}\def\dD{{\mathbb D}}\def\E{{\mathbb E}}\def\dF{{\mathbb F}}\def\H{{\mathbb H}}\def\dG{{\mathbb G}}
\def\I{{\mathbb I}}\def\J{{\mathbb J}}\def\K{{\mathbb K}}\def\dL{{\mathbb L}}\def\M{{\mathbb M}}
\def\N{{\mathbb N}}\def\dO{{\mathbb O}}\def\dP{{\mathbb P}}\def\Q{{\mathbb Q}}
\def\R{{\mathbb R}} \def\S{{\mathbb S}} \def\T{{\mathbb T}} \def\U{{\mathbb U}}
\def\V{{\mathbb V}}\def\W{{\mathbb W}} \def\X{{\mathbb X}}\def\Y{{\mathbb Y}}\def\Z{{\mathbb Z}}

%%%%% Arrows %%%%%
\def\la{\leftarrow}              \def\ra{\rightarrow}     \def\Ra{\Rightarrow}
\def\ua{\uparrow}                \def\da{\downarrow}
\def\lra{\leftrightarrow}        \def\Lra{\Leftrightarrow}
\newcommand{\abs}[1]{\lvert#1\rvert}
\newcommand{\br}[1]{\left(#1\right)}
\def\lan{\langle} \def\ran{\rangle}
%%%%% Typography %%%%%
\def\lt{\biggl}                  \def\rt{\biggr}
\def\ol{\overline}               \def\wt{\widetilde}
\def\no{\noindent}
%%%%% Math signs %%%%%
\let\ge\geqslant                 \let\le\leqslant
\def\lan{\langle}                \def\ran{\rangle}
\def\/{\over}                    \def\iy{\infty}
\def\sm{\setminus}               \def\es{\emptyset}
\def\ss{\subset}                 \def\ts{\times}
\def\pa{\partial}                \def\os{\oplus}
\def\om{\ominus}                 \def\ev{\equiv}
\def\iint{\int\!\!\!\int}        \def\iintt{\mathop{\int\!\!\int\!\!\dots\!\!\int}\limits}
\def\el2{\ell^{\,2}}             \def\1{1\!\!1}
\def\sh{\sharp}
\def\wh{\widehat}
\def\bs{\backslash}
\def\na{\nabla}
\def\ti{\tilde}
%%%%% Math operations %%%%%
\def\sh{\mathop{\mathrm{sh}}\nolimits}
\def\all{\mathop{\mathrm{all}}\nolimits}
\def\Area{\mathop{\mathrm{Area}}\nolimits}
\def\arg{\mathop{\mathrm{arg}}\nolimits}
\def\const{\mathop{\mathrm{const}}\nolimits}
\def\det{\mathop{\mathrm{det}}\nolimits}
\def\diag{\mathop{\mathrm{diag}}\nolimits}
\def\diam{\mathop{\mathrm{diam}}\nolimits}
\def\dim{\mathop{\mathrm{dim}}\nolimits}
\def\dist{\mathop{\mathrm{dist}}\nolimits}
\def\Im{\mathop{\mathrm{Im}}\nolimits}
\def\Iso{\mathop{\mathrm{Iso}}\nolimits}
\def\Ker{\mathop{\mathrm{Ker}}\nolimits}
\def\Lip{\mathop{\mathrm{Lip}}\nolimits}
\def\rank{\mathop{\mathrm{rank}}\limits}
\def\Ran{\mathop{\mathrm{Ran}}\nolimits}
\def\Re{\mathop{\mathrm{Re}}\nolimits}
\def\Res{\mathop{\mathrm{Res}}\nolimits}
\def\res{\mathop{\mathrm{res}}\limits}
\def\sign{\mathop{\mathrm{sign}}\nolimits}
\def\span{\mathop{\mathrm{span}}\nolimits}
\def\supp{\mathop{\mathrm{supp}}\nolimits}
\def\Tr{\mathop{\mathrm{Tr}}\nolimits}
\def\BBox{\hspace{1mm}\vrule height6pt width5.5pt depth0pt \hspace{6pt}}
\def\where{\mathop{\mathrm{where}}\nolimits}
\def\as{\mathop{\mathrm{as}}\nolimits}
%%%%%%%%%%%%% specialities %%%%%%%%%%%%%%
\newcommand\nh[2]{\widehat{#1}\vphantom{#1}^{(#2)}}
%{{\mathop{#1}\limits^\wedge}\vphantom{#1}^{(#2)}}
\def\dia{\diamond}
\def\Oplus{\bigoplus\nolimits}
%%%%%%%%%%% End of definitions %%%%%%%%%%
%%%%% OLD OLD OLD
\def\qqq{\qquad}
\def\qq{\quad}
\let\ge\geqslant
\let\le\leqslant
\let\geq\geqslant
\let\leq\leqslant
\newcommand{\ca}{\begin{cases}}
\newcommand{\ac}{\end{cases}}
\newcommand{\ma}{\begin{pmatrix}}
\newcommand{\am}{\end{pmatrix}}
\renewcommand{\[}{\begin{equation}}
\renewcommand{\]}{\end{equation}}
\def\eq{\begin{equation}}
\def\qe{\end{equation}}
\def\[{\begin{equation}}
\def\bu{\bullet}
\def\ced{\centerdot}
\def\tes{\textstyle}
\newcommand{\hb}{\hbar}
\newcommand{\Bk}{\Bbbk}

%\bigskip
%\medskip
%\smallskip

\title[{Resonances for vector-valued Jacobi
operators on half-lattice}] {Resonances for vector-valued Jacobi
operators on half-lattice}

\date{\today}
\author[Evgeny Korotyaev]{Evgeny Korotyaev }
\address{ Saint-Petersburg
State University, Universitetskaya nab. 7/9, St. Petersburg, 199034,
Russia and HSE University, 3A Kantemirovskaya ulitsa, St.
Petersburg, 194100, Russia, \ korotyaev@gmail.com, \
e.korotyaev@spbu.ru}

\begin{abstract}
\no  We study resonances for Jacobi operators on the half lattice
with matrix valued coefficient and  finitely supported
perturbations.  We describe a forbidden domain, the geometry of
resonances and their asymptotics when the main coefficient of the
perturbation (determining its length of
support)  goes to zero. Moreover, we show that\\
1) Any sequence of points on the complex plane can be resonances for
some Jacobi operators. In particular, the
multiplicity of a resonance can be any number.\\
2) The Jost determinant coincides with the Fredholm determinant up
to the constant.\\
3) The S-matrix on the a.c spectrum determines the perturbation
uniquely.\\
5) The value of  the Jost matrix at  any finite sequence of points
on the a.c spectrum   determine the Jacobi matrix uniquely. The
length of this sequence is equals to the upper point of the support
perturbation.

\end{abstract}
\maketitle

\begin{quotation}
\begin{center}
{\bf Table of Contents}
\end{center}

\vskip 6pt

{\footnotesize

1. Introduction and main results \hfill \pageref{Sec1}\ \ \ \ \

2. Preliminaries  \hfill \pageref{Sec2}\ \ \ \ \

3. Properties of the Jost solutions \hfill \pageref{Sec3}\ \ \ \ \

4. Proof of main theorems  \hfill \pageref{Sec4}\ \ \ \ \

5. Fredholm determinants \hfill \pageref{Sec5}\ \ \ \ \

6. Appendix: Inverse problems for finite Jacobi operators  \hfill
\pageref{Sec6}\ \ \ \ \

 }
\end{quotation}

\section {\lb{Sec1} Introduction and main results}
\setcounter{equation}{0}

\subsection{Introduction}

We consider Jacobi operators $H$ with matrix-valued coefficients
acting on $y=(y_x)_1^\iy\in \mH:= \ell^2(\N, \C^d), d\ge 1$ and
given by
\[
\lb{dH} (Hy)_x=a_{x-1}y_{x-1}+a_{x}y_{x+1}+b_xy_x, \qqq \ x\in
\N=\{1,2,3,...\},
\]
where formally $y_0=0$. We assume that $a_x>0, b_x=b^*$, $x\in\N$
are sequences of complex $d\ts d$ matrices. We discuss resonances
for finitely supported sequences $a_x-\1_d, b_x$, where $\1=\1_d$ is
the identity $d\ts d$ matrix.  There are a lot of results about
resonances for the scalar  case $d=1$, see \cite{BNW05},
\cite{DZ19}, \cite{IK10}--\cite{IK12} \cite{K11}, \cite{KL23},
\cite{MNSW12} and the references therein.
 There are results about direct and inverse scattering
theory for this class of operators with $d>1$, see for example,
 Aptekarev--Nikishin \cite{AN83},
Geronimo \cite{G82}, Nikishin \cite{Ni84, Ni87}  and Serebryakov
\cite{Se80, Se86} and recent papers \cite{BS21}, \cite{BS22}.
Different types of Jacobi operators were discussed in \cite{CGR05},
\cite{GKM02}, \cite{DS08}, \cite{KK08}, \cite{KK09}, \cite{Kn16},
\cite{Kn12}, \cite{Kn11}, \cite{S07}.

We consider also non-selfadjoint Jacobi operators $H_c=H_o+V_c $
acting on $ \mH$, where the unperturbed  operator $H_o$ and the
perturbation $V_c$ have the form
\[
\lb{HV}
\begin{aligned}
(H_oy)_x=y_{x-1}+y_{x+1},\qq  (V_cy)_x=(a_{x-1}-\1_d)y_{x-1}
 +b_xy_x+(s_x-\1_d)y_{x+1},\qq y_x\in\C^d,
\end{aligned}
\]
where the operator $H_o$ is the discrete Laplacian on $\N$.  We
assume that $a_x, s_x, b_x$, $x\in\N$ are sequences of complex $d\ts
d$ matrices. We assume that $V_c\in \J_{q}^c$  defined below. We
define the space $\M=\M_d$ of all $d\ts d$ complex matrices and its
subsets
\[
\begin{aligned}
\lb{dd} \M^s=\M_d^s=\{ b\in \M_d: b=b^*\},
 \qqq \M^o=\M_d^o=\{b\in \M_d: \det b\ne
0\}.
\end{aligned}
\]
 We study resonances for following classes of perturbations.

\no {\bf Definition V.} {\it  Introduce the class $\J_{q}^c, q\in
\{2p-1,2p\}, p\in \N$ of all operators $V_c$ given by \er{HV}, where
sequences $(b_x,a_{x}-\1_d,s_{1}-\1_d, x\in \N)$ are finitely
supported and satisfy
\[
 \lb{qs}
\ca  a_{x}, s_{x}\in \M_d^o, \ b_x\in \M_d, \ x\ge 1,
\\
a_x=s_x=\1_d, \  b_x = 0,\qq \forall \ x>p\ac, \qq
 \ca a_p=s_p=\1_d, b_p\in \M_d^o\ \ & {\rm if}\ q= 2p-1
 \\
    \1-s_pa_p \in \M_d^o,  \ \ & {\rm if}  \ q=2p   \ac ,
\]
and its subset of self-adjoint perturbations $V\in \J_{q}$ given by}
 \[
\lb{sa}  \J_{q}=\{V\in \J_q^c:   a_x=s_x>0, \qq b_x=b_x^*, \qq
\forall \ x\ge 1\}.
\]

Our main goal is to study  self-adjoint operators $H=H_o+V$, where
$V$ is given by
\[
\lb{Jo}  (Vy)_x=(a_{x-1}-\1_d)y_{x-1}+b_xy_x+
(a_x^*-\1_d)y_{x+1},\qq y_x\in \C^d.
\]
We consider resonances for two natural classes:\\
$\bu$ all matrices $a_x>0$,\\
$\bu$ all matrices $a_x$ are upper triangular with positive elements
on the diagonal.

The last case corresponds precisely Jacobi matrices with $2d + 1$
non-vanishing diagonals with the extreme diagonals strictly
positive. Without loss of the generality we assume that $V\in
\J_{q}$, since there are exact relations between these two classes,
see Sect. 6. The class $\J_{q}^c$ corresponds to non-self-adjoint
operators $H_c$ and the class $\J_{q}$ corresponds to self-adjoint
operators $H$.

A non-self-adjoint operator $H_c$ has an essential spectrum $[-2,2]$
plus a finite number $\gn_\bu$ of eigenvalues $\l_j$ (with algebraic
multiplicity $d_j$) on the set $\C\sm [-2,2]$. Recall the main
properties of the self-adjoint operator $H$ from \cite{AN83}. It has
purely absolutely continuous spectrum $[-2,2]$ plus a finite number
$\gn$ of eigenvalues $\l_j$ (with multiplicity $d_j\le d$) on the
set $\R\sm [-2,2]$. These eigenvalues of $H$ are denoted by
\[
\lb{eig1} \l_{-1_-}< \dots<  \l_{-\gn_-}<-2<2<\l_{\gn_+}<\dots<
\l_{1},\qqq \gn=\gn_++\gn_-,
\]
for some $\gn_\pm\ge 0$. Thus the spectrum of $H$ has the form (see
Fig. 1)
$$
\lb{sJ} \s(H)=\s_{ac}(H)\cup \s_d(H),\qqq \s_{ac}(H)=[-2,2],\qqq
\s_d(H)\ss \R\sm [-2,2].
$$
\setlength{\unitlength}{1.0mm}
\begin{figure}[h]
\centering
\unitlength 1.0mm % = 2.845pt
\begin{picture}(135,25)
\put(5,10){\line(10,0){100.00}}
\put(25,9){$\ts$} \put(25,13){$\l_{-1}$}
\put(51,9){$\ts$} \put(35,12.5){$\l_{-2}$}
\put(35,9){$\ts$}\put(51,12.5){$\l_{-3}$}
\put(56,6){-2} \put(58,10){\line(0,1){1.6}}
\put(58,11.6){\line(1,0){22.00}} \put(80,10){\line(0,1){1.6}}
\put(65,15){$\s_{ac}(H)$}
\put(78,6){2}
\put(85,9){$\ts$} \put(92,9){$\ts$} \put(85,12.5){$\l_{2}$}
\put(92,12.5){$\l_{1}$}
\end{picture}
%\vspace{-10mm}
\caption{\footnotesize The spectrum of $H$} \label{fig}
\end{figure}
 The operators $H$ and $ H_c$ with a perturbation $V\in
\J_{q}$ and $V_c\in \J_{q}^c$, with $q\in\{2p-1, 2p\}$, represent
 semi-infinite matrices of the form
%\[\label{J0}
$$
H=\left(\begin{array}{ccccccccc} b_1 & a_1 & 0 & 0 &  ...  &  0 & 0
& 0 \cr
 a_1 & b_2 & a_2 & 0 &  ... &  0& 0 & 0  \cr
 0 & a_2 & b_3 & a_3 &  ... & 0 & 0 & 0 \cr
 ... & ... & ... & ... & ... & ... & ... & ...  \cr
 0 & ...  & 0 & a_{p-1} & b_{p} & a_{p}& 0& 0 \cr
 0 & ...  & 0 & 0 & a_{p} & 0 & \1 & 0 \cr
 0 & ...  & 0 & 0 & 0 & \1 &0 & \1 \cr
  ... & ... & ... & ... & ... & ... & ... & ...
\end{array}\right),\
H_c=\left(\begin{array}{ccccccccc} b_1 & s_1 & 0 & 0 &  ...   & 0 &
0 & 0\cr
 a_1 & b_2 & s_2 & 0 &  ... &  0& 0 & 0  \cr
 0 & a_2 & b_3 & s_3 &  ... & 0 & 0 & 0 \cr
 ... & ... & ... & ... & ... & ... & ... & ... \cr
 0 & ...  & 0 & a_{p-1} & b_{p} & s_{p}& 0& 0 \cr
 0 & ...  & 0 & 0 & a_{p} & 0 & \1 & 0 \cr
 0 & ...  & 0 & 0 & 0 & \1 &0 & \1 \cr
  ... & ... & ... & ... & ... & ... & ... & ...
\end{array}\right).
$$
Define a standard spectral variable $k= k(\l)$ by
\[
\lb{lz0}
\textstyle \l=\l(k)=k+{1\/k}, \qqq k\in \dD, \;\; \l\in \L:=\C\sm
[-2,2],
\]
where $\;\dD_r=\{k\in \C:|k|<r\}, r>0$ is the disc, and $\dD =
\dD_1$. Here $\l(k)$ is a conformal mapping from $\dD$ onto the cut
domain $\L$ and its inverse function $k(\l)$ is
\[
\lb{kL} \textstyle
 k(\l)={\l\/2}-\sqrt{{\l^2\/4}-1},\qqq
k(\l)=\frac{1+o(1)}{\l}\;\;\; \text{ as } \l\to\iy .
\]
For an operator $H_c$ we introduce a Jost solution
$f(k)=(f_x(k))_0^\iy, f_x(k)\in \M_d$ to the equation
\[
\lb{jej}
\begin{aligned}
& a_{x-1}f_{x-1}+s_xf_{x+1}+b_xf_x=(k+k^{-1}) f_x,\qqq \ |k|>0,\;\;
x\in \Z_+=\{0,1,2,3,...\},
\\
& \where \qqq f_x(k)=k^{x} \1_d\ \qq \ \forall  x>p\qq {\rm and }
\qq q=\{2p-1,2p\}.
\end{aligned}
\]
For convenience, in order to define $f_0(k)$ we set $a_0:= \1_d$.
For the operator $H_c$ we define the $d\ts d$ Jost  matrix by
$\p=f_0(k)$, which is our basic function. The Jost matrix $\p(k)$ is
a $\M_d$-valued polynomial in $k\in \C$ of order $q$. Define  the
Jost determinant $\gf(k)=\det \p(k)$. The function $\gf(k)$ is a
polynomial of order $qd$ with zeros $r_1,....,r_{qd}$ labeled by
$$
0<|r_1|\le |r_2|\le \dots \le |r_{qd}|.
$$
Each zero $r_j\in \dD$ is an eigenvalue of $H_c$.  The function
$\gf$ can also have a finite number of zeros in $\C\sm \ol\dD$ which
are called resonances, plus possible zeros at $\pm 1$. The points
$\pm 1$ are  called virtual or resonances. In the self-adjoint case
$V\in \J_q$  the function $\gf$  has $\gn_\bu=\gn_++\gn_-\ge 0$ real
zeros in $\dD$ given by
\[
\lb{eg1}
\begin{aligned}
&\textstyle k_j =k(\l_j),\qq \l_j=k_j+{1\/k_j},\qqq j\in \gn_{\bu}:=
\{-\gn_-,\dots, -1, 1, \dots, \gn_+\},
\\
& -1<k_{-\gn_-}<...<k_{-1}<0<k_{1}<...<k_{\gn_+}<1.
 \end{aligned}
\]
Each $k_j$ has a multiplicity $d_j\le d$. The resonances of $H$ can
also be defined as the poles of  its resolvent meromorphically
continued in $\C$, see Proposition \ref{Trr1}.

For the operator $H_c$ given by \er{HV}, we define matrices
\[
\begin{aligned}
 \lb{m1x}
& \gt_x=a_x^{-1},\qq \cT_x=\gt_{1}\gt_{2}\cdot \cdot \cdot\gt_{x},
\qq c_x=1-s_xa_x,\ \ \qq \cC_0= \cT_{p-1}c_p\gt_{p},
\\
& \cC_0^o=\cT_{p-1}b_p\gt_{p} +\cT_{p-1}^1c_p\gt_{p}, \qqq \cC_0^1=
\cT_{p-1}^1b_p+\cT_{p-1} c_{p-1}, \qq B_x=\sum_1^x b_j,
\\
& \cT_{x}^1=b_1\cT_{x}+\gt_{1}b_2\gt_{2}...\gt_{x}+
\gt_{1}\gt_{2}b_3\gt_{3}...\gt_{x}+...+ \cT_{x-1}b_{x}\gt_{x},\qq
x\in \N_p=\{1,...,p\}.
\end{aligned}
\]
Our first results are devoted to the operator $H_c$ with a finitely
supported perturbation, including a selfadjoint case, when
$c_p=\1-a_p^2$. Define the ring of all zeros
$$
\mD_{ze}=\{r_-\le |r_j|\le r_+\}, \qq r_+=\max |r_j|,\qq r_-=\min
|r_j|.
$$
We describe the basic properties of Jost matrices  of the operator
$H_c$ and the forbidden domain for resonances defined by
$\mD_f=\{|k|>r_+\}$.

\begin{theorem}
\label{T1}
Let $V_c\in \J_q^c$ for some $q\ge1$. Then

\no i) The Jost matrix $\p=f_0(k)$ has asymptotics
\[
\lb{Ja1}
\begin{aligned}
\p(k)=\ca \cC_0k^{q}-\cC_0^ok^{q-1}+O(k^{q-2})\ \  & if \qq q=2p
\\
-\cT_{p-1}b_pk^q+\cC_0^1k^{q-1}+O(k^{q-2})\ \ & if \qq q=2p-1 \ac\qq
\as \qqq k\to \iy,
\end{aligned}
\]
and
\[
\lb{Ja2}
\begin{aligned}
& \p(k)=\cT_p-\cT_p^1 k+O(k^2)
\\
& \gf(k)= (1-k\Tr B_p+O(k^2))\det \cT_p,
\end{aligned}
\qqq  \as \qqq k\to 0.
\]
 \no ii) The Jost determinant $\gf=\det \p$ is a polynomial and has the form
\[
\lb{f02-4}
 \gf(k)=C_V\prod_1^{qd}(k-r_n),\qq C_V=(-1)^{dq}\gh \det\cT_p,\qq
\gh=\ca \det c_p \ {\rm if} \ q=2p\\
         \det b_p \ {\rm if} \ q=2p-1\ac,
\]
\[
\lb{f02-6}  \prod_1^{qd}r_n={1\/\gh},\qq
\sum_1^{qd}\frac{1}{r_n}=\sum_{x=1}^{p}\Tr b_x.
 \qq
\]
\no iii) If  the Jost determinant $\gf$ does not have zeros in $\dD$
for $q>2$, then
\[
\lb{f02-6x} 1\le r_-\le  |\gh|^{-{1\/qd}}\le r_+\le |\gh|^{-{1}}.
\]
\no iv) If $|\gh|>1$, then  $H_c$  has $\gm$ eigenvalues, counted
with multiplicities, and $r_-^\gm<1/|\gh|$.

\end{theorem}

 Define the resolvent $R_o(\l)=(H_o-\l)^{-1}$ and the   Fredholm
determinant $ D_c(k)$ by
\[
\lb{dD} D_c(k)=\det (I+V_cR_o(\l)), \qqq \tes \l=k+{1\/k},\ \  k\in
\dD.
\]
This determinant is well defined since $V_c$ is finitely supported.
The zeros of $D_c(k)$ in the disc $\dD$ are the zeros of $\gf$. We
will to show that $D_c=C\gf$ for some constant $C\ne0$. It is well
known that each $k_j\in \dD$ is a zero of $D_c$ and has a
multiplicity $d_j\le d$. The S-matrix for $H$ is defined by
\[
\lb{dsm} S(k)=\p^{-1}(k)\p(\ol k),  \ \ \ \ k\in \S:=\{|z| = 1\}.
\]
The function $S$ has a meromorphic continuation from the circle $\S$
to the whole complex plane. Then the poles of $S(k)$ are exactly the
resonances of $H$.

\begin{theorem}
\label{T1x} Let $V_c\in \J_q^c$ for some $q\ge1$. Then
 the Fredholm determinant $D_c(k), k\in
\dD$ defined by \er{dD} satisfies
\[
\lb{JfD} \gf(k)=D_c(k)\det\cT_p \qqq \forall \ \ \ k\in \dD,
\]
and then for all $k\in \C$. Moreover, if $V\in \J_q$, then the
function $\gf$ is real on the real line and the resonances are
symmetric with respect to the real line.
\end{theorem}

\no {\bf Remarks.}   The proof of the identity $\gf=D\det\cT_p$ for
$d=1$ is sufficiently simple. In the matrix  case it is not enough
and we need to use additional arguments, see Theorem \ref{Ta3}.

Now we consider the resonances for  the case $a_p=\ve a_p'$ where
the coupling constant $\ve\to 0$ while all other coefficients are
fixed, and $a_p', b_p\in \M_d^o$. Thus we discuss the case when the
support of the perturbation $V$ becomes smaller.

\begin{theorem}
\label{T2} Let $V\in \J_{q}, q\in \N$ satisfy
\[
\lb{asr1}
\begin{aligned}
  a_p=\ve a_p',\qq a_p',b_p\in \M_d^o\ \  {\rm if} \qq  q=2p
  \ \ \and \ \ b_p=\ve b_p',\qq b_p', c_{p-1}\in \M_d^o, \ \ {\rm if} \qq
q=2p-1,
\end{aligned}
\]
and  $\ve\to 0$. Let $\p(k,\ve)$ be the corresponding Jost matrix.
Let $k_j^o, j\in \N_m$ be zeros of $\gf(k,0)=\det\p(k,0)$ with the
multiplicity $d_j^o$ for some $N$. Then

\no i) The function $\p(k,\ve)$ converges to $\p(k,0)$ uniformly in
each disk as $ \ve\to 0$.

\no ii)  For any $\ve >0$ small enough there exist $d_j^o$
zeros(counted with multiplicity) of $\gf(k,\ve)=\det\p(k,\ve)$ in
each disc $\dD_{r_\ve}(k_j^o):=\{|k-k_j^o|<r_\ve\}$ for some
$r_\ve>0$ small enough.

\no iii) Moreover, $\gf(k,\ve)\ne 0$ in the domain $\dD_{R_\ve}\sm
\cup \dD_{r_\ve}(k_j^o)$ and  $\gf(k,\ve)$ has $d$ zeros
$k_{n}(\ve) \in \{|k|>R_\ve\}$ (counted with multiplicity) for some
$R_\ve>1$ big enough and they have asymptotics
\[
\lb{are1} k_{n}(\ve)=[\t_{n}+O(\ve^{{1\/d}})]/\ve, \qq n=1,...,d,
\]
where $\t_n, n\in \N_d$ are zeros of $\det(A+\t B), \t\in \C$, and
$A=2a'_p, B=b_p$ if $q=2p$ and $A=-b'_p, B=c_{p-1}$ if $q=2p-1$,
these zeros are real and $\ne 0$.

\end{theorem}

\no{\bf Remark.} Similar results for $d=1$ were obtained in
\cite{MW07}, when all $a_n-1, b_n\to 0$. Our case for $d=1$ was
discussed in \cite{KL23}.

\subsection{Inverse problems}

We discuss inverse problems for selfadjoint operators $H$. Firstly
we consider possible positions of eigenvalues and resonances.

\begin{theorem}\lb{T3}

i)  Let finite sequence of numbers $k_j\in (-1,0)\cup (0,1), j\in
N_m$ with the multiplicity $d_j\le d$  for some $m\ge 0$  be given.
Then there exists a finite perturbation $V\in \J_q$ for some $q\ge
1$ such that each $k_j, j\in N_m$ is an eigenvalue (with the
multiplicity $d_j$) of the operator $H$.

\no ii) Let a sequence  $r_n\in \C\sm \dD, n=1,..., qd$ be given for
some $q\ge 1$ and a number of the points $r_n=1$ (or $r_n=-1$) have
 multiplicity $\le d$. Then there exists a finite perturbation
$V\in \J_q$ such that each point $r_n$  is a resonance for the
operator $H$.

\no iii) The multiplicity of a resonance can be any number for some
specific perturbation $V\in \J_q$.

\end{theorem}

\no{\bf Remark.} The proof of the theorem is based on results from
\cite{KL23}, see Section 4.

We formulate our second result about inverse  problems  for Jacobi
operators  with finitely supported perturbations on the half
lattice. We show that the operator $H$  is uniquely determined by
the Jost  matrix at a finite number of points.

\begin{theorem}
\lb{T4}

Let an operator  $V\in \J_q$ for some $q\in \N$. Then

\no i)  The S-matrix  $S(k)$ for all $k\in \S$ determine $V\in \J_q,
q\in \N$ uniquely.

\no ii) Let $\z_j\in \C,\ j\in\N_{q+1}$ be a sequence of distinct
numbers. Then the sequence of the Jost  matrix $\p(\z_j),
j\in\N_{q+1}$ determine $V$ uniquely.

\end{theorem}

\no {\bf Remark.} The proof of Theorem \ref{T4} i) is based on
results from \cite{K19}, see Section 6.

\subsection{Historical review}
There are a lot papers devoted to resonances for various operators
see articles \cite{F97, H99, K04, S00, Z87} concerning to
one-dimensional Schr\"odinger operator with compactly supported
potentials and the book \cite{DZ19} and the references therein. In
particular, the inverse resonances problem (uniqueness,
reconstruction and characterization) were solved by Korotyaev for a
Schr\"odinger operator with a compactly supported real-valued
potential on the half-line \cite{K04} and on the real line
\cite{K05}. See also Zworski \cite{Z01} and Brown-Knowles-Weikard
\cite{BKW03} concerning the uniqueness problem. The resonance
stability was considered in papers \cite{K04a, MSW10, B12}.
Resonances are discussed for various perturbations of
one-dimensional Schr\"odinger operator. For example, Sch\"odinger
operator with periodic plus a compactly supported potential were
considered by Firsova \cite{F84}, Korotyaev \cite{K11} and
Korotyaev-Schmidt \cite{KS12}. Resonances for Schr\"odinger operator
with a step-like potential were discussed by Christiansen
\cite{C06}. Stark operator with compactly supported potentials was
studied in \cite{FH21, K17}. Resonances were also considered for
three and fourth order differential operator with compactly
supported coefficients on the line, see \cite{K19, BK19}. Note that
the inverse resonance scattering for Dirac operators was discussed
in \cite{IK14,IK14a,    KM21, KM23,M22}.

Now we discuss the case of Schr\"odinger operators on the cubic
lattice $\Z^d, d\ge 2$ for self-adjoint potentials. For
Schr\"odinger operators with decaying potentials on the lattice
$\Z^d$, Boutet de Monvel and Sahbani \cite{BS99} used Mourre's
method to prove completeness of the wave operators, absence of
singular continuous spectrum and local finiteness of eigenvalues
away from threshold energies. Isozaki and Korotyaev \cite{IK12}
studied the direct and the inverse scattering problem as well as
trace formulas. Korotyaev and Moller \cite{KM19} discussed the
spectral theory for potentials $V\in \ell^p, p>1$. An upper bound on
the number of discrete eigenvalues in terms of potentials was given
by Korotyaev and Sloushch \cite{KS20}, Rozenblum and Solomyak
\cite{RoS09}. Different types of trace formulas for Schr\"odinger
operators on discrete periodic graphs are discussed in  \cite{K24},
\cite{KL18}, \cite{KS22}, \cite{KS23}.

Finally,  we briefly comment on background literature for  Jacobi
operators. There are a lot of articles about resonances and inverse
resonance theory for scalar Jacobi operators with finitely supported
perturbation is well understood, see \cite{BNW05}, \cite{IK10},
\cite{IK11}, \cite{IK12A}, \cite{IK12B}, \cite{K11}, \cite{KL23},
\cite{Kn16}, \cite{MNSW12}, \cite{MW07},  and see references
therein. Inverse inverse resonance theory for scalar Jacobi
operators with finitely supported perturbation are discussed in
\cite{BNW05}, \cite{K11}, \cite{KL23} and the stability problem are
are discussed in  \cite{MW07}, \cite{KL23}. The forbidden domain for
resonances and the localization of complex and real resonances are
considered in \cite{KL23}.

The corresponding theory for Jacobi matrices with matrix-valued
coefficients is still largely a wide open field. The scattering
theory for Jacobi matrices with matrix-valued coefficients is
considered in \cite{AN83}, \cite{G82}, \cite{Ni84}, \cite{Ni87},
\cite{Se85}, \cite{Se86}. The inverse scattering problems for such
operators are discussed in \cite{Se80}, \cite{Se86}. The scattering
matrix for Schr\"odinger operators with matrix-valued potential on
the lattice are discussed in \cite{BS21}, \cite{BS22}.  Some
particular results were  obtained  in \cite{GKM02}, \cite{CGR05},
\cite{KK09}, \cite{KK08}, \cite{K21a}, \cite{DS08} and see
references therein.

\section {\lb{Sec2} Preliminaries}
\setcounter{equation}{0}

\subsection{Determinants and Traces}
Let $|\cdot|$ be the operator norm in the space $\C^d$.
 For a set $Y\subset\Z$ (below  $Y=\N,\; \Z_+=\N\cup \{0\}$ or
$\Z$) define the spaces $\ell^2(Y)$ and $\ell^{\iy }(Y)$ as the
spaces of sequences $(h_x)_{x\in Y}$ equipped with the norms given
by $ \|h\|_2^2 =\sum_{x\in Y}|h_x|^2\geq 0$ and $\|h\|_{\iy } =
\sup_{x\in Y}|h_x| $ respectively.

Let us recall some well-known facts from \cite{GK69}. Define the
space $\B$ of bounded operators acting on the Hilbert space $\cK$
equipped with the norm $ \|\cdot \|$. Let $\B_2(\cK)$ and
$\B_1(\cK)$ be the Hilbert-Schmidt class and the trace  of compact
operators on the Hilbert space $\cK$
 equipped with the norms
$\|A\|_{\B_2}^2=\Tr (A^*A)<\infty $ and  $ \|\cdot \|_{\B_1}$
respectively.

\no $\bu$ Let $A, B\in \B$ and $AB, BA\in \B_1$. Then
\begin{equation}
\label{AB} {\rm Tr}\, AB={\rm Tr}\, BA,
\end{equation}
\begin{equation}
\label{1+AB} \det (I+ AB)=\det (I+BA).
\end{equation}

\no $\bu$ Let $A\in \cB_1$. Then
\begin{equation}
\label{DA1} |\det (I+ A)|\le e^{\|A\|_{\cB_1}}.
\end{equation}

\no $\bu$  Let a function $A(\cdot): \mD  \to \B_1$ be analytic in a
domain $\mD \subset {\C}$, and the operator $(I+A(k))^{-1} $ is
bounded  for any $k\in\mD$. Then the function $F(k)=\det (I+ A(k))$
is analytic in $\mD$ and  satisfies
\[
\lb{m1}
 F'(k)= F(k)\Tr (I+A(k))^{-1}A'(k)\qqq \forall \ k\in\mD.
\]
\no $\bu$ Let $A\in \B_1$. Then the series $\sum _{j=1}^\iy \Tr
{(-A)^j\/j}$ converges and
\[
\lb{m2}
\det (I+A)=\exp \rt[- \sum _{j=1}^\iy \Tr {(-A)^j\/j}\rt],
\]
and, in particular,  if $\|A\|_{\B_1}$ is small enough, then
\[
\lb{m3} \ln \det (I+A)=\Tr \Big [A- {A^2\/2} +
O(\|A\|_{\B_2}^3)\Big].
\]

\no $\bu$ Let $A^{adj}$ be the adjugate matrix of a $d\ts d$ matrix
$A$. If $A$ is an invertible matrix, then
\[
\lb{13} A^{-1}={1\/\det A} A^{adj}.
\]
\no $\bu$ For a $d\ts d$ matrix $A$ the determinants has a
decomposition
\[
\lb{14} \det(A-\t \1_d)=\sum_0^{d}\gT_j\t^{d-j},
\]
 where  $\gT_0=1$  and the functions $\gT_j$ are given by (see p.331-333 in \cite{RS78})
\[
\lb{15}  \gT_1=-T_1,\ \  \gT_2=-{1\/2}(T_2+T_1\gT_1),\ \ .. ... ,
\gT_j=-{1\/m}\sum_0^{j-1}T_{j-s}\gT_s,.., \qq T_j=\Tr A^j, \qq j\in
\N.
\]

\subsection{Jost functions}

Recall that for operators $H_c$ given by \er{HV} we have defined
matrices
\[ \lb{f01-3}
\begin{aligned}
& \gt_x=a_x^{-1},\qq \cT_x=\gt_{1}\gt_{2}\cdot \cdot \cdot\gt_{x},
\qq \cT_{x,p}=\gt_{x}\gt_{2}\cdot \cdot \cdot\gt_{p}, \qq
c_x=1-s_xa_x,\ \ \qq \cC_0= \cT_{p-1}c_p\gt_{p},
\\
& \cC_0^o=\cT_{p-1}b_p\gt_{p} +\cT_{p-1}^1c_p\gt_{p}, \qqq \cC_0^1=
\cT_{p-1}^1b_p+\cT_{p-1} c_{p-1}, \qq B_x=\sum_1^x b_j,
\\
& \cT_{x}^1=b_1\cT_{x}+\gt_{1}b_2\gt_{2}...\gt_{x}+
\gt_{1}\gt_{2}b_3\gt_{3}...\gt_{x}+...+ \cT_{x-1}b_{x}\gt_{x},\qq
x\in \N_p=\{1,...,p\}.
\end{aligned}
\]
We discuss asymptotics of Jost solutions $f=(f_x)$ for
non-self-adjoint operators, in general.

\begin{lemma}\lb{TL1}
 Let an operator $V_c\in \J_q^c$ for some $q\ge1$.
 Then each function $f_x(k)$ is a
  matrix-valued polynomial in $k\in \C$  having the form
\[
\lb{f01-1} f_{p}(k)=\gt_p k^{p},\qqq
 f_{j}(k)=k^j\mF_{j}(k),  \qq j=0,1..., p-1,
\]
where $\mF_j(k)$ is a polynomial of a degree $q-2j$, and $f_{0}(k),
f_{1}(k)$ have asymptotics  as $k\to \iy$:
\[
\lb{f01-2}
\begin{aligned}
&  f_{0}(k)=\cC_0k^{q}-\cC_0^ok^{q-1}+O(k^{q-2})
\\
& f_{1}(k)=\cC_0k^{q-1}-(\cC_0^o-b_1\cC_0)k^{q-2}+O(k^{q-3})
\\
& f_{x}(k)=\gt_{x}...\gt_{p-1}c_p\gt_pk^{q-x}+O(k^{q-x-1}) ,\qq
x=2,...p-1,
 \end{aligned},
 \qqq q=2p,
\]
and
\[
\lb{f01-4xx}
\begin{aligned}
& f_{0}(k)=-\cT_{p-1}b_pk^q+\cC_0^1k^{q-1}+O(k^{q-2})
\\
&
f_{1}(k)=-\cT_{p-1}b_pk^{q-1}+(\cC_0^1-b_1\cT_{p-1}b_p)k^{q-2}+O(k^{q-2})
\\
& f_{x}(k)=-\gt_{x}...\gt_{p-1}b_pk^{q-x}+O(k^{q-x-1}),\qq
x=2,...p-1,
\end{aligned}
,\qqq q=2p-1.
\]
\end{lemma}

{\no\bf Proof.} Let $q=2p$. Equations
$a_{x}f_{x}=(\l-b_{x+1})f_{x+1}-s_{x+1}f_{x+2}$ and $f_x=k^x\1,\;
x>p$, give:

\no $\bu$ if $x=p$, then
\[
\lb{asp0} a_p f_p=(\l-b_{p+1})f_{p+1}-f_{p+2}
=(k+k^{-1})k^{p+1}\1-k^{p+2}\1=k^{p}\1\qqq \Rightarrow \qq
f_p=k^{p}\gt_p,
\]
\no $\bu$ if $x=p-1$, then
$$
\begin{aligned}
 & a_{p-1}f_{p-1}=(\l-b_{p})f_{p}-s_pf_{p+1}
=((k+k^{-1}-b_{p})k^{p}-s_pa_pk^{p+1}) \gt_p
=(k^{p+1}c_p-b_pk^{p}+k^{p-1})\gt_p,
\\
&f_{p-1}=\gt_{p-1}(k^{p+1}c_p-b_pk^{p}+k^{p-1})\gt_p=
\gt_{p-1}k^{p-1}(k^{2}c_p-b_pk+1)\gt_p,
\end{aligned}
$$
\no $\bu$ if $x=p-2$, then
$$
\begin{aligned}
& a_{p-2}f_{p-2}={(\l-b_{p-1})f_{p-1}-s_{p-1}f_{p}}
=\Big[(k+k^{-1}-b_{p-1})\gt_{p-1}(k^{p+1}c_p-b_pk^{p}+
k^{p-1})-s_{p-1}k^{p} \Big]\gt_p
\\
&= \Big[(k^{p+2}+k^{p})\gt_{p-1} c_p-k^{p+1}b_{p-1}\gt_{p-1}c_p-
(k^{p+1}+k^{p-1})\gt_{p-1}b_{p}+b_{p-1}\gt_{p-1}b_pk^{p}
\\
&+\gt_{p-1}(k^{p}+k^{p-2})
-b_{p-1}\gt_{p-1}k^{p-1}-s_{p-1}k^{p}\Big]\gt_{p}.
\end{aligned}
$$
 Thus we have
$$
\begin{aligned}
& f_{p-2}=\gt_{p-2}\Big[k^{p+2}\gt_{p-1}c_p-
k^{p+1}\Big(b_{p-1}\gt_{p-1}c_p +\gt_{p-1}b_{p} \Big)
\\
& +k^{p}\Big(\gt_{p-1}c_p + b_{p-1}\gt_{p-1}b_p+\gt_{p-1}-s_{p-1}
\Big)-k^{p-1}\Big(\gt_{p-1}b_p+  b_{p-1}\gt_{p-1}\Big)+
\gt_{p-1}k^{p-2}  \Big]\gt_{p},
\end{aligned}
$$
and defining  $\G_{p-j}=\gt_{p-j}\cdots\gt_{p-1}$ we get
$$
\begin{aligned}
 f_{p-2}
=k^{p+2}\G_{p-2}c_p\gt_{p}-k^{p+1}\cC_{p-2}^o+ O(k^{p}),\qqq
\cC_{p-2}^o=\gt_{p-2}\gt_{p-1}b_{p}\gt_{p}+\gt_{p-2}b_{p-1}\gt_{p-1}c_p
\gt_{p}.
\end{aligned}
$$
\no $\bu$ If $x=p-3$, then
$$
\begin{aligned}
&
a_{p-3}f_{p-3}=(\l-b_{p-2})f_{p-2}-s_{p-2}f_{p-1}=(k-b_{p-2})f_{p-2}+O(k^{p+1})
\\
& =(k-b_{p-2})\Big[k^{p+2}\G_{p-2}c_p\gt_p-
k^{p+1}\cC_{p-2}^o\Big]+O(k^{p+1})
\\
= & k^{p+3}\G_{p-2}c_p\gt_p-
k^{p+2}b_{p-2}\gt_{p-2}\gt_{p-1}c_p\gt_p- k^{p+2}\cC_{p-2}^o
+O(k^{p+1}).
\end{aligned}
$$
Thus we have
$$
\begin{aligned}
f_{p-3}=k^{p+3}\G_{p-3}c_p\gt_p-\cC_{p-3}^o k^{p+2}+O(k^{p+1}),
\\
\cC_{p-3}^o=\gt_{p-3}\gt_{p-2}\gt_{p-1}b_{p}\gt_p+
\gt_{p-3}[b_{p-2}\gt_{p-2}\gt_{p-1}+\gt_{p-2}b_{p-1}\gt_{p-1}]c_p\gt_p.
\end{aligned}
$$
Repeating this procedure we obtain \er{f01-1} and  \er{f01-2} for
$f_0$. We show \er{f01-2} for $f_1$. Above considerations give as
$k\to \iy$:
$$
f_{0}(k)=\cC_0k^{q}-\cC_0^ok^{q-1}+O(k^{q-2}),\qq
f_{1}(k)=\cC_ok^{q-1}-Ak^{q-2}+O(k^{q-3}), \qq f_{2}(k)=O(k^{q-3})
$$
for some matrix $A$. Then substituting these asymptotics  into the
equation for $f_0, f_1, f_2$ we have
\[
\begin{aligned}
\tes f_0=(k-b_1)f_1+({1\/k}f_1-s_1 f_2)=(k-b_1)f_1+O(k^{q-2}),
\end{aligned}
\]
which yields $A=\cC_0^o-b_1\cC_0$.

\no   Let $k=2p-1$ and $a_p=s_p=\1$. Then we have $\gt_p=\1$ and
$f_p=k^p\1$. Let $x=p-1$, then
$$
\begin{aligned}
 & a_{p-1}f_{p-1}=(\l-b_{p})f_{p}-f_{p+1}
=(k+k^{-1}-b_{p})k^{p}-k^{p+1} =-b_pk^{p}+k^{p-1},
\\
&f_{p-1}=\gt_{p-1}(-b_pk^{p}+k^{p-1}),
\end{aligned}
$$
Note that if $q=2p$ then  $f_{p-1}\sim k^{p+1}$ and but if $q=2p-1$
then $f_{p-1}\sim k^{p}$.

\no $\bu$  if $x=p-2$, then
$$
\begin{aligned}
& a_{p-2}f_{p-2}={(\l-b_{p-1})f_{p-1}-s_{p-1}f_{p}}
=(k-b_{p-1})\gt_{p-1}(-b_pk^{p}+ k^{p-1})-s_{p-1}k^{p}+O(k^{p-1})
\\
&=-k^{p+1}\gt_{p-1}b_{p}+k^{p}(b_{p-1}\gt_{p-1}b_p +\gt_{p-1}
-s_{p-1})+O(k^{p-1}).
\end{aligned}
$$
 Thus due to $\gt_{x}
-s_{x}=\gt_{x}c_{x}$ we have
$$
\begin{aligned}
 f_{p-2}
=-\G_{p-2}b_{p}k^{p+1}+\cC_{p-2}^1 k^{p}+ O(k^{p-1}),\qqq
\cC_{p-2}^1=
\gt_{p-2}\Big[b_{p-1}\gt_{p-1}b_p+\gt_{p-1}c_{p-1}\Big].
\end{aligned}
$$
$\bu$ If $x=p-3$, then due to $f_{p-1}=O(k^{p})$ we have
$$
\begin{aligned}
& a_{p-3}f_{p-3}=(\l-b_{p-2})f_{p-2}+O(k^{p})
%\\ &
=(k-b_{p-2})\Big[-k^{p+1}\G_{p-2}b_{p}+\cC_{p-2}^1
k^{p}\Big]+O(k^{p})
\\
&=- k^{p+2}\G_{p-2}b_{p} +
k^{p+1}\Big[b_{p-2}\G_{p-2}b_{p}+\cC_{p-2}^1 \Big]+O(k^{p}).
\end{aligned}
$$
Thus we have
$$
\begin{aligned}
f_{p-3}=- k^{p+2}\G_{p-3}b_{p} +\cC_{p-3}^1 k^{p+1}+O(k^{p}), \qq
\cC_{p-3}^1=\gt_{p-3}b_{p-2}\G_{p-2}b_{p}+\gt_{p-3}\cC_{p-2}^,
\end{aligned}
$$
Repeating this procedure we obtain \er{f01-1} and \er{f01-4xx} for
$f_0$.

 We show \er{f01-4xx} for $f_1$. Above considerations give as $k\to \iy$:
$$
f_{0}(k)=-\a k^{q}+\a^ok^{q-1}+O(k^{q-2}),\qq f_{1}(k)=-\a
k^{q-1}+Bk^{q-2}+O(k^{q-3}), \qq f_{2}(k)=O(k^{q-3})
$$
for some matrix $B$ and $\a=\cT_{p-1}b_p, \ \a^o=\cC_0^1$. Then the
substitution of these asymptotics into the equation for $f_0, f_1,
f_2$ implies
$$
\begin{aligned}
\tes f_0=(k-b_1)f_1+({1\/k}f_1-s_1 f_2)=(k-b_1)f_1+O(k^{q-2}),
\end{aligned}
$$
which yields $B=\a^o-b_1\a$. \BBox

We discuss the asymptotics of Jost solutions  for non-selfadjoint
operators  $H_c$ as $k\to 0$ and properties of the Jost function
$\gf(k)=\det f_{0}(k)$.

\begin{lemma}\lb{TL2}
 Let  $V\in \J_q^c$ for some $q\ge1$. Then
the polynomial $f_{0}(k)$ and Jost function $\gf(k)=\det f_{0}(k)$
have the following asymptotics as $k\to0$:
\[
\lb{f02-2}
\begin{aligned}
& f_{0}(k)=\cT_p-\cT_p^1 k+O(k^2),
\\
& f_{1}(k)=\cT_pk-(\cT_p^1-b_1\cT_p)k^2+O(k^3),
\\
& f_{2}(k)=\cT_{2,p}k^2+O(k^3),
\\
& \gf(k)=(1-k\Tr B_p+O(k^2))\det \cT_p,
\end{aligned}
\]
where $\cT_p=\gt_{1}\gt_{2}\cdot \cdot \cdot\gt_{p}$ and $
\cT_p^1=b_1\cT_p+\gt_{1}b_2\gt_{2}...\gt_{p}+...+
\cT_{p-1}b_{p}\gt_{p}$ and
\[
\lb{f02-3x} \Tr B_x=\sum_1^x \Tr
b_j=\Tr \cT_x^{-1}\cT_x^1, \qq x\in \N.
\]
\no  ii) Moreover, the determinant $\gf=\det f_0$ and its zeros
$r_j$ satisfy
\[
\lb{z1}
 \gf(k)=C_V\prod_1^{qd}(k-r_n),\qq C_V=(-1)^{dq}\gh \det\cT_p,\qq
\gh=\ca \det c_p \ {\rm if} \ q=2p\\
         \det b_p \ {\rm if} \ q=2p-1\ac,
\]
\[
\lb{f01-4q}
\begin{aligned}
 \gf(k)=C_V (k^q-k^{q-1}A+O(k^{q-2})),\ \
 A=\Tr \ca c_p^{-1}b_p+\sum_1^{p-1}b_j \
& if\ q=2p
\\
b_p^{-1}c_{p-1}+\sum_2^{p-1}b_j, \ & if \  q=2p-1\ac,
\end{aligned}
\]
\[
\lb{z2}
 \prod_1^{qd}r_j={1\/\gh},\qqq \sum_1^{qd}{1\/r_j}=\Tr B_p\qqq
 \sum_1^{qd}r_j=A.
\]

\end{lemma}

{\no\bf Proof.} i) Using the equation
$a_{x}f_{x}=(\l-b_{x+1})f_{x+1}-s_{x+1}f_{x+2}$ and $f_x=k^x,\;
x>p$, we obtain $f_p=k^{p}\gt_p$ and the following cases:

\no $\bu$ if $x=p-1$, then
$$
\begin{aligned}
 & a_{p-1}f_{p-1}=(\l-b_{p})f_{p}-s_pf_{p+1}
=((k+k^{-1}-b_{p})k^{p}-s_pa_pk^{p+1}) \gt_p
=(k^{p+1}c_p-b_pk^{p}+k^{p-1})\gt_p,
\\
&f_{p-1}=\gt_{p-1}(k^{p-1}-b_pk^{p}+O(k^{p+1}))\gt_p;
\end{aligned}
$$
$\bu$ if $x=p-2$, then
$$
\begin{aligned}
& a_{p-2}f_{p-2}={(\l-b_{p-1})f_{p-1}-s_{p-1}f_{p}}
\\
&=\Big[(k+k^{-1}-b_{p-1})\gt_{p-1}(k^{p-1}-b_pk^{p}
+O(k^{p+1}))+O(k^{p}) \Big]\gt_p
\\
&= \Big[\gt_{p-1}k^{p-2}-k^{p-1}b_{p-1}\gt_{p-1}-
k^{p-1}\gt_{p-1}b_{p}+O(k^{p})\Big]\gt_{p}.
\end{aligned}
$$
Thus we have
$$
\begin{aligned}
 f_{p-2}=k^{p-2}\gt_{p-2}\gt_{p-1}\gt_{p}- \cK_{p-2}\gt_{p}k^{p-1}+
 O(k^{p}),
\qq \cK_{p-2}=\gt_{p-2} b_{p-1}\gt_{p-1}+\gt_{p-2}\gt_{p-1}b_{p}.
\end{aligned}
$$
$\bu$ If $x=p-3$, then
$$
\begin{aligned}
 f_{p-3}=\gt_{p-3}\Big\{(\l-b_{p-2})f_{p-2}-s_{p-2}f_{p-1}\Big\}
\\
=\gt_{p-3}\Big\{(k+k^{-1}-b_{p-2})\Big[k^{p-2}\gt_{p-2}\gt_{p-1}-
 \cK_{p-2}k^{p-1}+O(k^{p})\Big]+O(k^{p-1})      \Big\}\gt_p
\\
= k^{p-3}\gt_{p-3}\gt_{p-2}\gt_{p-1}\gt_{p}-\cK_{p-3}\gt_{p}k^{p-2}
+O(k^{p-1}),
\end{aligned}
$$
where
$$
\cK_{p-3}= \gt_{p-3}\cK_{p-2}+\gt_{p-3}b_{p-2}\gt_{p-2}\gt_{p-1}=
\gt_{p-3}b_{p-2}\gt_{p-2}\gt_{p-1}+
\gt_{p-3}\gt_{p-2} b_{p-1}\gt_{p-1}+\gt_{p-3}\gt_{p-2}\gt_{p-1}b_{p}.
$$
Repeating this procedure we obtain asymptotics of $f_0$ in
\er{f02-2} and using \er{m3} we obtain   asymptotics of $\det f_0$
in \er{f02-2}, since the   identity \er{f02-3x} holds true:
$$
\Tr \cT_x^1\cT_x^{-1} =\Tr (
b_1\cT_x+\gt_{1}b_2\gt_{2}...\gt_{x}+...+ \cT_{x-1}b_{x}\gt_{x}
)\cT_x^{-1}=  \sum_1^x \Tr b_j=\Tr B_x.
$$
 We show \er{f02-2} for $f_1$. Above considerations give as $k\to 0$:
$$
f_{0}(k)=\cT_p-\cT_p^1 k+O(k^2),\qq f_{1}(k)=\cT_p
k-Bk^{2}+O(k^{3}), \qq f_{2}(k)=\cT_{2,p}k^2+O(k^3)
$$
for some matrix $B$. Then the substitution of these asymptotics into
the equation for $f_0, f_1, f_2$ implies
\[
\begin{aligned}
\tes f_0=(k+{1\/k}-b_1)f_1-s_1
f_2=({1\/k}-b_1)f_1+O(k^{2})=\cT_p-(B+b_1\cT_p)k+O(k^{2}),
\end{aligned}
\]
which yields \er{f02-2} for $f_1$.

\no ii) Asymptotics \er{f01-2}, \er{f01-4xx}  give \er{z1}.\\
Let $m=p-1$. If $q=2p$, then from \er{f01-2}, \er{z1} we obtain
$A=\Tr \cC_0^{-1}\cC_0^o$ and
$$
\begin{aligned}
A=\Tr [\cT_{m}c_p\gt_{p}]^{-1} \rt(\cT_{m}b_p\gt_{p}
+\cT_{m}^1c_p\gt_{p}\rt)= \Tr
[c_p^{-1}b_p+\cT_{m}^{-1}\cT_{m}^1]=\Tr
[c_p^{-1}b_p+\sum_1^{p-1}b_j]
\end{aligned}
$$
Let $q=2p-1$. Then from  \er{f01-4xx}, \er{z1} we obtain $A=\Tr
(\cT_{p-1}b_p)^{-1}(\cC_0^1-b_1\cT_{p-1}b_p)$ and
$$
\begin{aligned}
A=\Tr \big((\cT_{m}b_p)^{-1}[\cT_{m}^1b_p+\cT_{m} c_{m}]-b_1\big)
=\Tr [\cT_{m}^{-1}\cT_{m}^1+b_p^{-1}c_{m}-b_1]=\Tr[
b_p^{-1}c_{p-1}+\sum_2^{p-1}b_j].
\end{aligned}
$$
Collecting last identities we obtain \er{f01-4q}.

We show \er{z2}. Let  $P(k)=\prod_1^{qd}(k-r_j)$. From \er{f02-2}
and \er{z1} we obtain
$$
\gf(k)=C_V P(k),\qq  \gf(0)=\det\cT_p=C_V\prod_1^{qd}(-r_j),
$$
where $C_V=(-1)^{dq}\gh \det\cT_p$,  which yields
$\prod_1^{qd}r_j={1/\gh}$ and $P(0)={(-1)^{dq}/\gh}$.

Recall Vieta's formulas for our polynomial
$P(k)=\prod_1^{dq}(k-r_j)$:
\[
\begin{aligned}
P(k)=k^{dq}-k^{{dq}-1}C_1+.....-C_2k+P(0), \qq C_1=\sum_1^{dq}r_j,
\qqq C_2=\sum _1^{dq} {P(0)\/r_j}.
\end{aligned}
\]
Then using $\det \cT_p=C_VP(0)$ and    Vieta's formulas and
asymptotics \er{f02-2} we obtain
$$
\begin{aligned}
& \gf(k)=C_V P(k)=C_V(P(0)-kC_2+....),
\\
& \gf(k)=(\det \cT_p)( 1-k\Tr B_p+O(k^2))  =C_V(P(0)-k P(0)\Tr
B_p+O(k^2)),
\end{aligned}
$$
which yields $C_2=\sum_1^{qd}{P(0)\/r_j}=P(0)\Tr B_p$ and thus
$\sum_1^{qd}{1\/r_n}=\Tr B_p$, since $P(0)\ne 0$.

From \er{f01-4q} and Vieta's formulas we obtain
$$
\begin{aligned}
\gf(k)=C_V  (k^q-k^{q-1}A+O(k^{q-2})=C_V  (k^q-k^{q-1}C_1+O(k^{q-2})
\end{aligned}
$$
which yields $C_1=\sum_1^{qd}r_j=A$.
 \BBox

For Schr\"odiger operators $q=2p-1$, all $a_x=\gt_x=\cT_x\1$ and our
matrices have the form
\[ \lb{f01-3s}
\begin{aligned}
 c_x=0=\cC_0,\ \cC_0^o=b_p, \ \ \cT_{x}^1=B_x,\qq
\cC_0^1=B_{p-1}b_p, \ \ \cC_0^1-b_1\cT_{p-1}b_p= \sum_2^{p-1}b_jb_p.
\end{aligned}
\]
Now we collect properties of the Jost functions for Schr\"odinger
operators, i.e., when all $a_x=\1$. From Lemma \ref{TL1} and Lemma
\ref{TL2} we obtain the following .

\begin{corollary}{\bf (Schr\"odinger
operators)} \lb{Tz1} Let $q=2p-1$ and all $a_x=s_x=\1$. Then
\[
\lb{c1}
\begin{aligned}
& f_{0}(k)=-b_pk^q+B_{p-1}b_pk^{q-1}+O(k^{q-2})
\\
& f_{1}(k)=-b_pk^{q-1}+(B_{p-1}-b_1)b_pk^{q-2}+O(k^{q-2})
\end{aligned},
\qq as\qq k\to \iy,
\]
\[
\lb{c2}
\begin{aligned}
f_{0}(k)=\1-B_p k+O(k^2)
\end{aligned}
\qq as\qq k\to 0,
\]
where $B_x=\sum_1^xb_j$ and
\[ \lb{c3}
\begin{aligned}
\gf(k)=C_V\prod _{j=1}^q (k-r_j),\qqq C_V=(-1)^d \det (b_p),
\end{aligned}
\]
\[
\lb{c4}
\begin{aligned}
 \gf(k)=C_V (k^q-k^{q-1}B_{p-1}+O(k^{q-2})),\ \
   \qq \as \ \  k\to \iy
 \\
\gf(k)= (1-k\Tr B_p+O(k^2))\det \cT_p, \qq \as \ \ k\to 0.
\end{aligned}
\]

 \end{corollary}

\subsection {Regular solutions}

Define a polynomial
\[
\lb{pb1}
\b_x(\l)=\gs_x(\l-b_{x})\gs_{x-1}(\l-b_{x-1})....\gs_1(\l-b_{1}),
\]
where $\gs_x=s_{x}^{-1}$. The direct computation gives that the
polynomial $\b_x(\l)$ has asymptotics
\[
\lb{fi2}
\begin{aligned}
& \b_x(\l) =\cS_x\l^x-   \l^{x-1}\cS_x^o+O(\l^{x-2}),
\\
&  \cS_x=\gs_x\gs_{x-1}...\gs_1,\qq
\cS_x^o=\gs_xb_x\gs_{x-1}...\gs_1+\gs_x\gs_{x-1}b_{x-1}...\gs_1+.....
+\gs_x\gs_{x-1}...\gs_1b_1.
\end{aligned}
\]
Define the regular solutions $\vp, \vt$ of the equation  under the
initial conditions:
\[
\lb{fs}
\begin{aligned}
a_{x-1}y_{x-1}+s_xy_{x+1}+b_xy_x=\l y_x, \\ \vp_0=\vt_1 = 0,\qqq
\vp_1 =\vt_0= \1.
\end{aligned}
\]

\begin{lemma}
\lb{ldop}
Let $V\in \J_q^c$ for some $q\ge1$. Then the regular solutions
$\vp_{x+1}(\l), x\in \N_p$ satisfy
\[
\lb{fi1} \textstyle
\begin{aligned}
\vp_{x+1}(\l)
=\b_x(\l)+O(\l^{x-2})=\cS_x\l^x-\l^{x-1}\cS_x^o+O(\l^{x-2})
\end{aligned}
\]
as $ \l\to+\iy$, where  $\cS_x=\gs_x\gs_{x-1}...\gs_1$ and $\cS_x^o$
are defined by \er{fi2}.
\end{lemma}

\no{\bf Proof.}    Using  the equation
$s_{x-1}\vp_{x}=(\l-b_{x-1})\vp_{x-1}-a_{x-2}\vp_{x-2}$  we obtain
for $x=2,3$:
$$
 \textstyle
\vp_2=\gs_1[(\l-b_{1})\vp_{1}] =\gs_1(\l-b_{1}), \qq
\vp_{3}=\gs_2[(\l-b_{2})\vp_{2}-a_1] =\b_2(\l)-\gs_2a_1.
$$
since $\vp_0 = 0,\; \vp_1 = \1$.  If $x=4, 5$, then we have
$$
\begin{aligned}
\textstyle  \vp_{4}=\gs_3[(\l-b_{3})\vp_{3}-a_2\vp_{2}]
=\b_3(\l)+O(\l),
\\
\textstyle \vp_{5}=\gs_4[(\l-b_{4})\vp_{4}-a_3\vp_{3}]
=\b_4(\l)+O(\l^2).
\end{aligned}
$$
Repeating this procedure and using \er{fi2}, we obtain \er{fi1}.
\BBox

\subsection{\lb{Ex}Examples}

We discuss two simplest examples at $q=1$ and $q=2$.

\no {\bf Example 1, $q=1$. }  Consider the case $q=p=1, a_1=\1, \det
b_1\ne 0$. Consider the zeros $k_j, j\in \N_d$ of the Jost function
$\gf=\det (\1-b_1k)$. Then we have
\[
\p=\1-b_1k,\qqq k_j^{-1}, j\in\N_d \qq {\rm are\ eigenvalues\ of\ }
\ b_1.
\]
\no If $|k_j|>1$, then $k_j$ is a resonance of $H$.

\no If $|k_j|<1$, then $k_j$ is an eigenvalue. If $|k_j|=1$, then
$k_j$ is a virtual state.

\no If the eigenvalues of $b_1$  are $\pm 1$. Then the operator $H$
has only virtual states $\pm 1$, there are no resonances and
eigenvalues.

\no {\bf Example 2, $q=2$.}  Consider the case $q=2,p=1$.   In this
case $\det c_1\ne 0$, where $c_1=\1-a_1^2$.  Then from the proof of
Lemma \ref{TL1} at $x=p-1$ we obtain
\[
\lb{ex1} \p=(k^{2}c_1-b_1k+\1) a_{1}^{-1}, \qq \gf=\det \p=\det
(k^{2}c_1-b_1k+\1)\det a_{1}^{-1}.
\]
\no $\bu$  Let $b_1=0$ and let $\pm k_{j}, j\in \N_{d}$ be zeros of
the Jost function $\gf=\det (k^{2}c_1+\1)\det (a_{1})^{-1}$, since
the function    $\gf$ is even. Then each $k_j^{-2}$ is real
eigenvalues of $-c_1$ and $k_j^{-2}\in (-1,\iy)$. Then  the zeros of
$\gf$ satisfy

If $k_j^2\in (0,1)$,  then $\pm k_{ j}$ is an eigenvalue of $H$.

If $k_j^2\in (1,\iy)$,  then $\pm k_{j}$ is a resonance of $H$.

If $k_j^2\in (-\iy, -1)$,  then $\pm k_{j}\in i\R\sm [-1,1]$ is a
resonance of $H$.

If $k_j^2=1$,  then $\pm k_{j}\in \{\pm 1\}$ is a resonance of $H$.

\no $1_{\bu}$ If $b_1\ne 0$,  then we obtain
\[
\gf(k)= \det (k^{2}c_1-b_1k+\1)=k^{2d}(\1+O(1/k))\det c_1
\]
 Then the Jost function $\gf(k)$ has 2d zeros on the plane $\C$.
 Even in this simple case very difficult to determine eigenvalues
 and resonances of the operator $H$ for large $d$.

\no $2_{\bu}$ If $b_1=0$ and $a_1=\sqrt 2 \1_d$, then we have
$c_1=-\1$ and $\gf=\det [(\1-k^{2})\1] 2^{-d/2}$.  Thus all zeros of
$\gf$ are $\pm 1$ only. Then the operator $H$ does no have
resonances and eigenvalues, all zeros of $\gf$ are virtual states
$\pm 1$  only.

\section {\lb{Sec3} Properties of the Jost solutions}
\setcounter{equation}{0}

\subsection{Jost solutions and their Wronskian}
In this section we discuss   Jost functions for self-adjoint
perturbations $V\in \J_q$. Define two mapping for sequences of
matrix-valued polynomials $y(k)=(y_x(k))_0^{\iy}$ by
$$
\begin{aligned}
\wt y(k)=(\wt y_x({k}))_0^\iy,\qq \wt y_x({k}):=y_x^*(\ol{k}),\qq
\\
y^\ced(k)=(y_x^\ced({k}))_0^\iy,\qq y_x^\ced({k}):=y_x(\Bk),\qq
\Bk=1/k.
\end{aligned}
$$
For a Jacobi operator $H$ and for two sequences of matrix-valued
polynomials $y(k)=(y_x(k))_0^{\iy}$, $\e=(\e_x(k))_0^{\iy}$ we
introduce the Wronskian
\[
\lb{dw}
\begin{aligned}
& \{y,\e\}_x(k)=\wt y_x({k}) a_x \e_{x+1}(k)- \wt y_{x+1}({k})
a_x\e_{x}(k),\qqq x=0,1,2,...
\end{aligned}
\]
The Wronskian has the property
\[
\wt {\{y,\x\}}=-\{\x,y\}.
\]
If both $y$ and $\e$ are solutions of the Jacobi equation \er{jej},
then the Wronskian $\{y,\e\}_x(k)$ does not depend on $x$. Recall
that $a_0=\1$ by our convention. In particular, we obtain
\[
\label{wr1}
\begin{aligned}
\{\vt,\vt\}_x=0,\qq
 \{\vp,\vp\}_x=0,
\qqq
 \{\vt,\vp\}_x=\1,\qq x\ge 0,
\end{aligned}
\]
and
\[
\label{wr2}
\begin{aligned}
&  \{f,f\}_x=\wt f_xa_xf_{x+1}- \wt f_{x+1}a_xf_{x}=0,
\\
&   \{f^\ced,f\}_x=\wt f_x^\ced a_xf_{x+1}- \wt f_{x+1}^\ced
a_xf_{x}=k-\Bk,
\end{aligned}
\]
where $ f_x^\ced(k)=f_x(\Bk)$ and
\[
\lb{L4-5}
\begin{aligned}
& \{f,\vp\}_x=\wt f_xa_x\vp_{x+1}-\wt  f_{x+1}a_x\vp_{x}=\wt \p,
\\
& \{f^\ced,\vp\}_x=\wt f_x^\ced a_x\vp_{x+1}- \wt f_{x+1}^\ced
a_x\vp_{x}=\wt \p^\ced,
\end{aligned}
\]
where we omit the variables $\l$ and $k$ for shortness.

Define the Weyl-Titchmarsh functions
\[
\lb{df1}
\begin{aligned}
\textstyle \gF=\vt+\vp m_+, \qqq  m_+=f_1f_0^{-1}.
 \end{aligned}
\]

\begin{lemma}\lb{Tjf}
Let  a perturbation $V \in \J_{q},\; q\in \N$. Then $\gF, m_+$
satisfy
\[
f=\vt \p+\vp f_1=\gF f_0,
\]
\[
\lb{wt1}
\begin{aligned}
&  \{\gF,\gF\}=\{\gF,f\}=m_+-\wt m_+=0,
\\
&  \{\gF,\vp\}=\1.
\end{aligned}
\]
\end{lemma}
\no{\bf Proof.} Identities \er{L4-5} follow from the definitions of
solutions $f, \vp$. We have $f=\vt A+\vp B$ for some matrices $A,
B$. Then \er{wr1} imply
$$
\begin{aligned}
-\p=\{\vp,f\}=\{\vp,\vt\}A+\{\vp,\vp\}B=-A.
\end{aligned}
$$
Similar arguments imply $B=f_1$. From the properties of $\vp, \vt$
we obtain $\{\gF,\gF\}=m_+-\wt m_+$ and the identity \er{wr2} gives
\er{wt1}. This and the identity $\gF=ff_0^{-1}$ give $\{\gF,f\}=0$.
\BBox

Now we present the regular solutions $\vp_x$ in terms of the Jost
solutions.

\begin{lemma}\lb{TL4}
Let  a perturbation $V \in \J_{q},\; q\in \{2p-1, 2p\}$ and
$\Bk=1/k$. Then
 \[
 \lb{L4-1}\tes
 \vp=(k-\Bk)^{-1}(f\wt\p^\ced -f^\ced \wt\p),
 \]
 \[
 \lb{L4-2}\tes
\p\wt\p^\ced =\p^\ced \wt\p,
 \]
\[
\lb{L4-3} \tes
 \vp_{p}={1\/k-\Bk}a_p^{-1}(k^p\wt\p^\ced -\Bk^p\wt\p),
\]
 \[
\lb{L4-4} \tes
 \vp_{x}={1\/k-\Bk}(k^x\wt\p^\ced -\Bk^x\wt\p), \qq x>p,
\]
\end{lemma}

\no{\bf Proof.} We have $\vp_x=f_xC_1+f_x^\ced C_2$ for some
matrices $C_1,C_2$. Then \er{wr2} imply
$$
\begin{aligned}
\wt\p=\{f,\vp\}_x=\{f,f \}_xC_1+\{f,f^\ced \}_x C_2=\{f,f\}_xC_1
+\{f,f^\ced \}_xC_2=(\Bk-k)C_2,
\\
\end{aligned}
$$
which yields $C_2=-\wt\p(k)/(k-\Bk)$. Similar arguments imply
$C_1=\wt\p ^\ced /(k-\Bk)$. Identity \er{L4-1} gives
\er{L4-2}-\er{L4-4}.
 \BBox

 Now we present the regular solutions $\vt_x$ in terms of the Jost
solutions.

\begin{lemma}\lb{TLt}
Let  a perturbation $V \in \J_{q},\; q\in \{2p-1, 2p\}$ and
$\Bk=1/z$. Then
\[
 \lb{Lt1}\tes
 \vt ={1\/k-\Bk}(f\wt f_1-f^\ced \wt f_1^\ced ),
 \]
\[
\lb{Lt2} \tes
 \vt_{p}(\l) ={1\/k-\Bk}a_p^{-1}(k^p\wt f_1-\Bk^p \wt f_1^\ced ).
 \]
\end{lemma}

\no {\bf Proof.} We have $\vt(\l)=f C_1+f^\ced  C_2$, where $f^\ced
=f(1/k), k\in\C$ and $C_1, C_2$ are some matrices. The Wronskian
gives at $x=0$:
$$
\begin{aligned}
\{f,\vt\}_x=\wt f_x a_x\vt_{x+1}-\wt f_{x+1}a_x\vt_{x}=-\wt f_{1},
\\
\{f^\ced,\vt\}_x=\wt f_x^\ced a_x\vt_{x+1}-\wt f_{x+1}^\ced
a_x\vt_{x}=-\wt f_{1}^\ced,
\end{aligned}
$$
since $\vt^\ced =\vt$ and then \er{wr2} gives
$$
\begin{aligned}
-\wt f_{1}=\{f,\vt\}_x=\{f,f C_1+f^\ced C_2\}_x=\{f,f\}_x
C_1+\{f,f^\ced \}_x C_2=(\Bk-k)C_2,
\\
\end{aligned}
$$
which yields $C_2=\wt f_{1}/(k-\Bk)$. Similar arguments imply
$C_1=\wt f_{1}^\ced/(k-\Bk)$. From \er{Lt1} we have \er{Lt2}. \BBox

Now we write down the Jost function in terms of the regular
solutions.

\begin{lemma}\lb{Tw}
Let  a perturbation $V \in \J_{q},\; q\in \{2p-1, 2p\}$ and
$\l=k+{1\/k}$. Then
 \[
 \lb{w1}
 \wt\p(k) =k^p \vp_{p+1}(\l) - a_pk^{p+1}\vp_p(\l),
\]
 \[
 \lb{w2}
\tes
 f_1(k) = -k^p\vt_{p+1}(\l) + a_pk^{p+1}\vt_p(\l).
 \]
\end{lemma}

\no{\bf Proof.} From \er{L4-5} at $x=p$ we have \er{w1}. From Lemma
\ref{TLt}  we have \er{w2}. \BBox

\subsection {Spectral parameters}
We discuss spectral parameters of the operator $H$ associated its
with eigenvalues.

\no {\bf Definition of Spectral data.} {\it For each eigenvalue
$k_j\in (-1,1)$, $j\in n_\bu$ with the multiplicity $d_j\le d$, we
define the spectral data:

\no $\bu$ {\bf the subspaces}
\[
\label{cEDef}
\begin{aligned}
{E_j}=\Ker f_{o}(k_j) \ss\C^d,\qq  {E_j^\bu}=\Ker f_{o}^*(k_j)
\ss\C^d,\qq \dim E_j=\dim E_j^\bu=d_j\le d;
\end{aligned}
\]
$\bu$ {\bf the orthogonal projectors} $P_j:\C^d\to E_j$  and
$P_j^\bu:\C^d\to E_j^\bu$,
$$
P_j\C^d=E_j,\qqq P_j^\bu\C^d=E_j^\bu;
$$
\no $\bu$  {\bf the positive norming matrix} $P_jK_j P_j: E_j\to
E_j$, where the matrix $K_j$ is given by }
\[
\label{Kdef}
\begin{aligned}
K_j=\sum_{x\ge 1} f_x(k_j)^*f_x(k_j)>0,\qqq j\in n_{\bu}.
\end{aligned}
\]
We define also the matrix-valued functions of $k\in \C$:
\[
\label{uv} u_j(k)=P_j^\bu \p(k)P_j: \ E_j\to E_j^\bu,\qqq \
v_j(k)=P_j^\bu f_{1}(k)P_j: \ E_j\to E_j^\bu.
\]

\begin{lemma}
\lb{f02}

Let $V\in \J_q,\; q\in\{2p,2p-1\}$ and let $\{k_{j}, j\in \gn_\bu\}
$ be all zeros of $\det f_{o}(k)$ on $(-1,1)$ and let $
\Bk_j=k^{-1}_j$. Then for any $j\in n_{\bu}$ we have
\[
\lb{nm1}
  K_j=\sum_{x\ge 1}f_x^*(k_j)f_x(k_j)={1\/\dot \l(k_j)}\Big(\dot
f_{o}^*(k_j)f_1(k_j)-\dot f_1^*(k_j) f_{o}(k_j)\Big)  > 0,
\]
\[
\lb{nm2}
  K_j \ge C_j\1_d,\qqq C_j={k_j^{2p+2}\/1-k_j^2},
\]
\[
\lb{nm4} f_{1}(k_j)P_j\ss E_j^\bu,\qqq u_j(k_j)=f_1(k_j)P_j,
\]
\[
\lb{nm3}
  P_jK_jP_j={1\/\dot \l(k_j)}P_j\dot
\p^*(k_j)f_1(k_j)P_j ={1\/\dot \l(k_j)}\dot u_j^*(k_j) v_j(k_j)  \ge
C_jP_j.
\]
\end{lemma}
\no{\bf Proof.} Let $\dot u={\pa \/\pa z}u$ and $\l=k+k^{-1}, k\in
\R$. From \er{jej} we have
$$
\begin{aligned}
& a_{x-1}f_{x-1}+a_x\psi_{x+1}+b_xf_x=\l f_x,\qqq  x\in \Z_+,
\\
& a_{x-1}\dot f_{x-1}+a_x\dot f_{x+1}+b_x\dot f_x=\l \dot f_x+\dot
\l f_x,
\\
&    \dot f_{x-1}^*a_{x-1}+\dot f_{x+1}^*a_x+\dot f_x^*b_x=\l \dot
f_x^*+\dot \l f_x^* ,
\end{aligned}
$$
multiplying the first equation by $\dot f_x^*$ and the third
equation by $f_x$ and taking the difference
\[
\lb{jejz}
\begin{aligned}
& \dot f_x^* a_{x-1}f_{x-1}+\dot f_x^*a_xf_{x+1}+\dot f_x^*b_xf_x=\l
\dot f_x^*f_x,\qq x\in \Z_+,
\\
& \dot f_{x-1}^*a_{x-1}f_x+\dot f_{x+1}^*a_xf_x+\dot f_x^*b_xf_x=\l
\dot f_x^*f_x+\dot \l f_x^*f_x ,
\end{aligned}
\]
we obtain
$$
\begin{aligned}
\dot \l f_x^*f_x=\dot f_{x-1}^*a_{x-1}f_x+\dot f_{x+1}^*a_xf_x-\dot
f_x^* a_{x-1}f_{x-1}-\dot f_x^*a_x\psi_{x+1}.
\end{aligned}
$$
Summing we have
$$
\begin{aligned}
\dot \l \sum_{x\ge 1}f_x^*f_x=\sum_{x\ge 1}\rt([\dot
f_{x-1}^*a_{x-1}f_x-\dot f_x^*a_x\psi_{x+1} ]+[\dot
f_{x+1}^*a_xf_x-\dot f_x^* a_{x-1}f_{x-1}]  \rt)
\\
=\dot f_{0}^*a_{0}f_1-\dot f_1^* a_{0}f_{0}=\dot f_{0}^*f_1-\dot
f_1^* f_{0},
\end{aligned}
$$
which yields \er{nm1}. We show \er{nm1}. From \er{jej} we have
$$
\sum_{x\ge p+1}f_x^*(k_j)f_x(k_j)=k_j^{2p+2}\sum_{n\ge
0}k_j^{2n}\1_d ={k_j^{2p+2}\/1-k_j^2}\1_d,
$$
which yields \er{nm2} and \er{nm3}. We show  \er{nm4}.
From \er{wr2} we obtain
$$
\begin{aligned}
\p^*(k_j)f_{1}(k_j)P_j=f_{1}^*(k_j)f_{0}(k_j)P_j=0,
\end{aligned}
$$
this gives $f_{1}(k_j)P_j\ss E_j^\bu$, i.e., $P_j^\bu
 f_{1}(k_j)P_j=f_{1}(k_j)P_j$, which yields  $P_j\dot
f_{0}^*(k_j)f_1(k_j)P_j =\dot u_j^*(k_j) v_j(k_j)$. Substituting it
into \er{nm1} and using \er{nm2} we obtain \er{nm3}.
 \BBox

We describe the poles of the matrix $\p(k)^{-1}$.

 \begin{lemma}
\label{Tin} Let $V\in \J_q,\; q\in \N$ and let $k_{j}\in (-1,1),
j\in \gn_\bu $ be a zero of $\gf$. Then

\no i)  The matrix $\dot u_j(k_j)=P_j^\bu\dot{\p}(k_j)P_j$ and
$P_j^\bu{\p}_{1}(k_j)P_j$ from $E_j$ to $E_j^\bu$ are invertible.

\no ii) The following asymptotics holds true:
\[
\label{po} \p(k)^{-1} = \frac{P_j\dot
u_j(k_j)^{-1}P_j^\bu}{k-k_j}+O(1)\qq {as}\qq k\to k_j.
\]
\end{lemma}
\no {\bf Proof.} i) Using (\ref{nm2}) and (\ref{nm4}), we obtain
$$
 C_jP_j\le  P_jK_jP_j={1\/\dot \l(k_j)}\dot u_j^*(k_j) v_j(k_j).
$$
Since the left hand side is positive definite on $E_j$, then the
matrices $v_j(k_j)=P_j^\bu f_1(k_j)P_j$   and $\dot
u_j(k_j)=P_j^\bu\dot{\p}(k_j)$ acting from $E_j$ to $E_j^\bu$ are
invertible.

\no ii)  Let $t=k-k_j\to 0$. Below we omit $j$ for shortness. Let $
\p= \ma {\a} \ {\b}\\ {\g} \ {u}\am $, where
$$
\begin{aligned}
\a_t=(P_j^\bu)^\bot \p(k_j+t)P_j^\bot:&
E_j^\bot\to(E_j^\bu)^\bot,\qq &
\b_t=(P_j^\bu)^\bot \p(k_j+t)P_j:& E_j\to( E_j^\bu)^\bot, \cr
\g_t=P_j^\bu \p(k_j+t)P_j^\bot:& E_j^\bot\to E_j^\bu\ , &
u_t=P_j^\bu \p(k_j+t)P_j:& E_j\to E_j^\bu\ .
\end{aligned}
$$
Due to $\p(k_j)P_j=0$ and $P_j^\bu \p(k_j)=0$, we have
$$
\a_t=\a_o+O(t), \qqq \b_t= t \dot \b_o+O(t^2),\qq \g_t= t \dot
\g_o+O(t^2),\qq u_t= t \dot u_o+O(t^2),
$$
as $t\to0$, where
$$
\dot u_o=P_j^\bu\dot{\p}(k_j)P_j=\dot u(k_j),\qqq
\a_o=(P_j^\bu)^\bot \p(k_j)P_j^\bot=\p(k_j)P_j^\bot:
E_j^\bot\to(E_j^\bu)^\bot.
$$
The operator $\a_o$ is invertible, since $\p h=0$ implies $h \in
E_j$. Due to i) the operator $\dot u_o$ is invertible.  Define
$w=u-\g \a^{-1}\b$ from the properties of $\a, \b, \g, u$ we deduce
that
$$
w_t=t(\dot u(k_j)+O(t)).
$$
Therefore, the standard formula for the inverse matrix
$$
\ma \a & \b \\ \g & u\am ^{-1}=
\ma {\a^{-1}+\a^{-1}\b w^{-1} \g \a^{-1}} & -\a^{-1}\b w^{-1} \\
 {-w^{-1}\g\a^{-1}} & w^{-1}\am,\qq w=u-\g\a^{-1}\b,
$$
\[
\det \ma \a & \b \\ \g & u\am=\det \a \det(u-\g\a^{-1}\b)= \det u
\det(\a-\b u^{-1}\g)
\]
gives
$$
\p(k)^{-1} = \ma {\a_o^{-1}+O(t)} & {O(1)} \\ {O(1)} & {t^{-1}\dot
u_o^{-1}+O(1)}\am \qq {\rm as} \qq t\to 0.
$$
In particular, this implies $\res_{k=k_j}(\p(k))^{-1}=P_j\dot
u_o^{-1}P_j^\bu$. \BBox

We describe properties of eigenvalues and the corresponding
projectors.

\begin{lemma}
\lb{Tpr1} Let $V\in \J_q,\; q\in\N$ and $k_{j}\in (-1,1), j\in
\gn_\bu $ be a zero of $\gf$ and $\Bk_j=1/k_j$.
 Then

\no i) The matrix $P_j^\bu \p(\Bk_j)$ satisfies
\[
\lb{p1x} \p^*(\Bk_j)P_j^\bu\ss E_j,\qqq
P_j\p^*(\Bk_j)P_j^\bu=u_j^*(\Bk_j)=\p^*(\Bk_j)P_j^\bu.
\]
\no ii) For the eigenvalue $k_j$ the corresponding solution has the
form
\[
 \lb{p2}\tes
 \vp(\l_j)P_j^\bu= \gc_jf(k_j)u_j(\Bk_j)^*,\qqq \gc_j={1\/k_j-\Bk_j},
 \]
\[
 \lb{p3}
\dot u_j(k_j)u_j(\Bk_j)^*=u_j(\Bk_j)\dot u_j(k_j)^*.
\]
\end{lemma}

\no {\bf Proof.} i)  From \er{L4-2} we have
$\p(k_j)\p^*(\Bk_j)P_j^\bu= \p(\Bk_j)\p^*(k_j)P_j^\bu=0$ which
yields $\p^*(\Bk_j)P_j^\bu\ss E_j$ and \er{p1x}.

\no ii) From \er{L4-1} we obtain
\[
 \lb{rs}\tes
 \vp(\l_j)=\gc_j(f(k_j)\p^*(\ol\Bk_j)-f(\Bk_j)\wt\p(k_j))=
 \gc_j(f(k_j)\p^*(\Bk_j)-f(\Bk_j)\p^*(k_j)).
 \]
Then from the identity \er{p1x}
$\p^*(k_j)P_j^\bu=0$ we obtain \er{p2}. From \er{L4-2}, \er{p1x} we
have
$$
\begin{aligned}
\tes P_j^\bu\dot\p(k_j)\p^*(\Bk_j)P_j^\bu=
P_j^\bu\p(\Bk_j)\dot\p^*(k_j)P_j^\bu= \dot
u_j(k_j)u_j(\Bk_j)^*=u_j(\Bk_j)\dot u_j(k_j)^*.
\end{aligned}
$$
\BBox

\section {\lb{Sec4} Proof of  inverse problems}
\setcounter{equation}{0}

\subsection {S-matrix}
We discuss properties of S-matrix $S(k)=\p^{-1}(k)\p(\ol k), \ \ \ \
k\in \S:=\{|z| = 1\}$. We determine asymptotics of S-matrix.

\begin{corollary}
\lb{TL3} Let $V\in \J_q$. Then S-matrix $S(k), k\in \S$ has a
meromorphic extension on the whole complex plane and satisfies
\[
\lb{sm} S(k)=\p^{-1}(k)\p(1/ k)=\wt\p(1/k)\wt\p(k)^{-1}, \qqq k\in
\C,
\]
\[
\lb{s1} S^*(k)= S(\ol k)= S(k)^{-1},\qqq k\in \S.
\]
Moreover, it has asymptotics:
\[
\lb{asS1z} S(k)=(-k)^{-q}\cC^{-1}\rt(\1+
S_1k^{-1}+O(k^{-2})\rt)\cT_p\qqq \as \qq k\to \iy,
\]
\[
\lb{asS2z} S(k)=(-k)^{-q}\cT_p^{-1}\rt(\1-k S_1+O(k^{2}\rt)\cC\qqq
\as \qq k\to 0,
\]

\[
\lb{cv1}
\begin{aligned}
\cC=\ca\cC_0 \\ \cT_{p-1}b_p\ac, \qq S_1=\ca
\cC_0^o\cC^{-1}-\cT_p^1\cT_p^{-1},& \qq if \qq q=2p
\\
\cC_0^1\cC^{-1}-\cT_p^1\cT_p^{-1},& \qq if \qq q=2p-1\ac
\end{aligned}
\]
where $\cT_p^1, \cT_p, \cC_0^o, \cC_0^1, \cC_0$ are given by
\er{m1x}.
 \end{corollary}
\no {\bf Proof.} Let $\Bk=1/k$. From the definition
$S(k)=\p^{-1}(k)\p(\ol k)$ and $\ol k=1/k$ on $\S$  and from
\er{L4-2} we obtain \er{sm}. From \er{sm} we get
 $$
S^*(k)=\p^{-1}(\ol k)\p( k) =S(\ol k),\qq
 S^*(k)S(k) =\p^{-1}(\ol k)\p( k)\p^{-1}(k)\p(\ol k)=\1.
 $$
\no $\bu$\ Let $q=2p$ and $\cC=\cC_0$. If $k\to \iy$, then from
$S(k)=\p(k)^{-1}\p(\Bk)$  and \er{f01-2}, \er{f02-2} we have :
$$
\begin{aligned}
S(k)={\p^{-1}(k)\p(\Bk)}=
k^{-q}\cC^{-1}[\1-\cC_0^o\cC^{-1}\Bk+O(\Bk^{2})]^{-1} [\cT_p-\Bk
\cT_p^1+O(\Bk^2)]
\\
=k^{-q}\cC^{-1}[(\1+\cC_0^o\cC^{-1}\Bk+O(\Bk^{2})] [\1-\Bk
\cT_p^1\cT_p^{-1}+O(\Bk^2)]\cT_p
\\
=k^{-q}\cC^{-1}\rt(\1+\Bk (\cC_0^o\cC^{-1}-\cT_p^1\cT_p^{-1})
+O(\Bk^{2})\rt)\cT_p.
\end{aligned}
$$
If $k\to 0$, then similar arguments give
$$
\begin{aligned}
S(k)=\p(k)^{-1}\p(\Bk)=
\Bk^{-q}\cT_p^{-1}[\1-\cT_p^{1}\cT_p^{-1}k+O(k^{2})]^{-1}
[\1-\cC_0^o\cC^{-1}k+O(k^{2})]\cC
\\
=\Bk^{2p}\cT_p^{-1}[\1+\cT_p^{1}\cT_p^{-1}k+O(k^{2})]
[\1-\cC_0^o\cC^{-1}k+O(k^2)] \cC
\\
=\Bk^{q}\cT_p^{-1}\rt(\1+k(\cT_p^{1}\cT_p^{-1}-\cC_0^o\cC^{-1})
+O(k^{2})\rt)\cC,
\end{aligned}
$$
\no $\bu$\  Let  $q=2p-1$ and $\cC=\cT_{p-1}b_p$. If  $k\to \iy$,
then from \er{f01-4xx}, \er{f02-2} we have
$$
\begin{aligned}
S(k)={\p(k)^{-1}\p(\Bk)}=-k^{-q}C^{-1}[\1-\cC_0^1C^{-1}\Bk+O(\Bk^{2})]
^{-1} [(\1-\Bk \cT_p^1\cT_p^{-1}+O(\Bk^2)]\cT_p
\\
=-k^{-q}C^{-1} [\1+\cC_0^1C^{-1}\Bk+O(\Bk^{2})] [\1-\Bk
\cT_p^1\cT_p^{-1}+O(\Bk^2)]\cT_p
\\
=-k^{-q}C^{-1}\rt(\1+\Bk ( \cC_0^1C^{-1} -\cT_p^1\cT_p^{-1})
+O(\Bk^{2})\rt)\cT_p.
\end{aligned}
$$
 If $k\to 0$, then similar arguments give
$$
\begin{aligned}
S(k)=\p^{-1}(k)\p(\Bk)= -k^{-q} \cT_p^{-1}[\1-k
\cT_p^1\cT_p^{-1}+O(k^2)]^{-1} [\1-\cC_0^1C^{-1}k+O(k^{2})]C
\\
=-k^{-q} \cT_p^{-1}[\1+k(\cT_p^1\cT_p^{-1}-\cC_0^1C^{-1})+O(k^2) ]C.
\end{aligned}
$$
\BBox

We rewrite the S-matrix in terms of the Weyl function
$M_j=\vp_{j}\vp_{j+1}^{-1}, j=p,p+1$ and conversly.

\begin{lemma}
\label{TSk}
 Let $V\in \J_q$ and $M_p(\l)=\vp_{p}(\l)\vp_{p+1}(\l)^{-1}$ and
 $\l(k)=k+{1\/k}, k\in\S\sm\{\pm 1\}$ and let $\vp_{j}(\l(k))\ne 0, j=p,
p+1,p+2$,. Then

\no i)  We have $\det(a_jM_j(\l(k))-\z)\ne 0$, where $j=p,p+1, $ and
$\z=k,\ol k$ and the identities hold true:
 \[
\lb{sm1} S(k)={1\/ k^{2p}}\cdot{\1-\Bk a_pM_p(\l(k))\/ \1-k
a_pM_p(\l(k))}={1\/ k^{2p+2}}\cdot{\1-\Bk M_{p+1}(\l(k))\/ \1-k
M_{p+1}(\l)}.
\]
ii)  The matrix $k^{2j}S(k)-\1, j=p,p+1,p+2$ is invertible and
\[
\lb{sm2} M_p(\l(k))=a_p^{-1}{k^{2p}S(k)-\1 \/k^{(2p+1)}S(k)-\Bk},
\qqq M_{p+1}(\l(k))={k^{2(p+1)}S(k)-\1\/k^{(2p+3)}S(k)-\Bk}.
\]
\end{lemma}
 \no {\bf Proof.}  i) Let $k\in \S$. From  \er{w1}  we have
 $$
 \begin{aligned}
& \wt\p(k) =k^p[\vp_{p+1}(\l) -
ka_p\vp_p(\l)]=k^p[(\1-ka_pM_p(\l)]\vp_{p+1}(\l),
\end{aligned}
 $$
 which implies  $\det(X_j(\l)-\z)\ne 0$.
Moreover, due to  the identity \er{sm} we obtain
$$
 \begin{aligned}
 S(k)=\wt\p(1/k)\wt\p(k)^{-1}={k}^{-2p}[\1- \Bk  a_pM_p(\l)]
[\1-k a_pM_p(\l)]^{-1}.
\end{aligned}
 $$
ii) From  \er{sm1}  we have $ k^{2p}S(k)[\1-k a_pM_p(\l)]=\1-\Bk a_p
M_p(\l) $, which yields \er{sm2}.
 \BBox

 We describe the poles of the matrix $S(k)$ associated with eigenvalues of
 $H$.

\begin{lemma}
\lb{Tpr1x} Let $k_j\in (-1,1)$ and $\Bk_j=1/k_j, j\in \gn_{\bu}$.
Then the matrix-valued function $S(k)$ has a simple pole at $k_j$
such that
\[
\lb{p1}
\begin{aligned}
& S(k)={X_j\/k-k_j}+O(1),
\\
& X_j=P_j\dot u_j(k_j)^{-1}P_j^\bu \p(\Bk_j)=P_j\dot u_j(k_j)^{-1}
u_j(\Bk_j)>0 \qq {\rm in }\qq E_j,
\end{aligned}
\]
where $\dot u_j(k_j)=P_j^\bu\dot{\p}(k_j)P_j:E_j\to E_j^\bu$ is
invertible for all $j\in n_\bu$.
\end{lemma}

\no {\bf Proof.}  From \er{po}, \er{p1x} we obtain as $k\to k_j$:
\[
\begin{aligned}
S(k)=\p^{-1}(k)\p(1/k)= \rt({P_j\dot
u_j(k_j)^{-1}P_j^\bu\/k-k_j}+O(1)\rt)\p(1/k)
\\
={X_j\/k-k_j}+O(1),\qqq X_j=P_j\dot u_j(k_j)^{-1}P_j^\bu \p(\Bk_j)=
P_j\dot u_j(k_j)^{-1}u_j(\Bk_j),
\end{aligned}
\]
where due to Lemma \ref{Tin} the mapping $\dot
u_j(k_j)=P_j^\bu\dot{\p}(k_j)P_j:E_j\to E_j^\bu$ is invertible for
all $j\in\N_m$. From \er{nm3}  we have
\[
\lb{1zz} C_jP_j\le   P_jK_jP_j={1\/\dot \l(k_j)}\dot u_j^*(k_j)
v_j(k_j) ={1\/\dot \l(k_j)}v_j^*(k_j) \dot u_j(k_j).
\]
From \er{wr2} and $u(k_j)P_j=0$ we obtain
$$
\begin{aligned}
&  \p^*(\Bk_j)f_{1}(k_j)- f_{1}^*(\Bk_j)\p(k_j)=k_j-\Bk_j,
\\
& P_j\p^*(\Bk_j)f_{1}(k_j)P_j=u_j^*(\Bk_j)v_j(k_j)
=v_j^*(k_j)u_j(\Bk_j)=(k_j-\Bk_j)P_j.
\end{aligned}
$$
Then $v_j^*(k_j)=(k_j-\Bk_j)P_ju_j(\Bk_j)^{-1}$ and substituting it
into \er{1zz} we obtain
$$
P_jK_jP_j={1\/\dot \l(k_j)}v_j^* (k_j)\dot u_j(k_j)={1\/\dot
\l(k_j)} (k_j-\Bk_j)P_ju_j(\Bk_j)^{-1}\dot
u_j(k_j)=k_jP_ju_j(\Bk_j)^{-1}\dot u_j(k_j).
$$
\BBox

\subsection {Proof of main theorems}

We begin to prove main theorems.

\no {\bf Proof of Theorem \ref{T1}.} i)  From \er{f01-2} and
\er{f02-2} we have \er{Ja1}, \er{Ja2} respectively.

\no ii) Let $\Bk=1/k$ and $q=2p$. We show \er{f02-4}. From \er{Ja1},
\er{m3} we obtain
 \[
 \begin{aligned}
\gf(k)=\det [\cC_0k^{q}+O(k^{q-1})]=k^{dq}\det \cC_0 \det
[\1+O(1/k)]=k^{dq}\det \cC_0  (\1+O(1/k)),
\end{aligned}
 \]
and $\det \cC_0=\det c_p \det\cT_p$. This yields $\gf(k)=\det c_p
\det\cT_p \prod_1^{qd}(k-r_n)$. The proof of \er{f02-4} for $q=2p-1$
is similar.

The asymptotics \er{Ja1} yields $\gf(0)=\det\cT_p \ne 0$, thus
$\prod_1^{qd}r_n={1\/\gh}$.
 From \er{Ja2} we have
$
{\gf'(0)\/\gf(0)}=-\sum_1^{k}\frac{1}{r_n}=-\Tr B_p,
$
which yields \er{f02-6}.

\no  iii) Recall that $\gh=\det c_p$ if $q=2p$ and $\gh=\det b_p$ if
$q=2p-1$. If the operator $H$ does not have eigenvalues for some
$q>2$, then all $|r_j|\ge 1$ (and some of them $>1$). From
\er{f02-4} we get $\prod_1^{qd}|r_j|=1/|\gh|$ and we deduce that
$|\gh|<1$ and $r_+<{1\/|\gh|}$ and $r_+^{dq}\le {1\/|\gh|}$, which
yields ${1\/|\gh|^{1\/qd}}\le r_+<{1\/|\gh|}$. Similar arguments
give $r_-^{qd}\le {1\/|\gh|}$. Collecting these estimates  we obtain
\er{f02-6x}.

\no iv) If $|\gh|>1$, then the identity
$\prod_1^{qd}|r_j|={1\/|\gh|}<1$ gives that some $|r_j|<1$ and the
operator $H$ has $\gm$ eigenvalues counted with multiplicities
eigenvalues. This yields   $r_-^\gm<1/|\gh|$.

Statement  v) is proved in Theorem \ref{Ta3}. \BBox

 Recall  results about the perturbation theory for polynomials, see
e.g.  p.5, from \cite{RS78}.

\begin{theorem}\lb{Tp1}

  Let $F(t,\ve)=t^d+v_1(\ve) t^{d-1}+..+v_d(\ve)$ be an d-th degree
polynomial in $t$ and whose coefficients are all analytic functions
of $\ve$. Suppose $t=t_o$ is a root of multiplicity $m$ of $F(t,0)$.
Then for $\ve$ near $\ve_o=0$, there are exactly $d_o$ roots
(counting multiplicity) of $F(t,\ve)$ near $t_o$ and these roots are
the branches of one or more multi-valued analytic functions with at
worst algebraic branch points at $\ve=\ve_o$. Explicitly, there are
positive integers $d_1,...,d_k$ with $d_1+...+d_k=d_o$ and
multi-valued analytic functions $t_1,...,t_k$ (not necessary
distinct) with convergent Puiseux series (Taylor series in
$\ve^{1\/d_o}$)
\[
\lb{Ps} t_s(\ve)=t_o+\sum_{j=1}^\iy \a_{s,j}\ve^{j\/d_s},
\]
so that the $d_o$ roots near $t_o$ are given by the $d_1$ values of
$t_1$, the $d_2$ values of $t_2$, etc.
\end{theorem}

{\no\bf Proof of Theorem \ref{T2}.} Consider $H$, where $V\in \J_q$
and $V$ satisfies \er{qs}, \er{asr1} for some $q=2p\ge2$. The proof
for the case $q=2p-1$ is similar. The condition \er{asr1} gives
$$
c_p=\1-a_p^2=-2\ve a_p'-\ve^2{a_p'}^2, \ \ \  \and \qq b_p, a_p'\in
\M_d^o.
$$
The corresponding Jost matrix we denote by $f_{0}(k,\ve)$. Let the
Jost polynomial $\gf(k,\ve)=\det \p(k,\ve)$ and $\p(k,0)=\det
\gf(k,0)$. Here
 $\gf(k,\ve)$ is a polynomial in $k$ of the order $qd$ and
 $\gf(k,0)$ is a polynomial in $k$ of the order $(q-1)d$.

i) Let $\gf(k,\ve)=\det \p(k,\ve)$ and $\gf(k,0)=\det \p(k,0)$. Here
 $\gf(k,\ve), \ve \ne 0$ is a polynomial in $k$ of the order $qd$ and
 $\gf(k,0)$ is a polynomial in $k$ of the order $(q-1)d$.
 Let $\ve \to 0$. Then the polynomial $\gf(k,\ve)$ converges to
the polynomial $\gf(k,0)$ uniformly in each disk $\dD_R, R>0$.

Let $k_j^o, j\in \N_M$ be zeros of $\gf(k,0)$ with the multiplicity
$d_j^o$ for some $M$. Consider the discs $\dD_{r_\ve}(k_j^o)$ for
some  $r_\ve>0$ small enough. Then due to by Rouch\'e's Theorem,
$\gf(k,\ve)$ has as many roots, counted with multiplicities, as
$\gf(k,0)$ in each disk $\dD_{r_ve}(k_j^o)$ and the remaining domain
$\dD_{R_\ve}\sm \cup \dD_{r_\ve}(k_j^o)$ for some $R_\ve>0$ big
enough.

\no ii) Consider the last $d$ zeros in the domain $\C\sm
\dD_{R_\ve}$ and show that they satisfy \er{are1}. In order to
compute these zeros we introduce new variable $z=1/k$ and the
corresponding functions
$$
\begin{aligned}
& \f(z,\ve)=z^{q}\p(1/z,\ve)=\cC_0-\cC_0^o z+z^2\f_1(z,\ve),
\\
& \cC_0=-\ve A-\ve^2 B,\qq      A= \cT_{p-1}2a_p'\gt_{p},\qq
B=\cT_{p-1}{a_p'}^2\gt_{p},
\\
& \cC_0^o=\a-\ve \b-\ve^2 \b_1,\qq  \a=\cT_{p-1}b_p\gt_{p}\in
\M_d^o,\qq \b= 2\cT_{p-1}^1 a_p'\gt_{p}, \qq \b_1= 2\cT_{p-1}^1
{a_p'}^2\gt_{p},
\end{aligned}
$$
where
\[
 |\f_1(z,\ve)|\le C_\f,
  \]
uniformly in $z, \ve\in \ol\dD$, for some constant $C_\f$. The large
zeros $k_j(\ve)\in \C\sm \dD_{R_\ve}, j\in \N_d$ map to zeros
$z_j(\ve)=1/k_j(\ve)\in \dD_{r_\ve}, j\in \N_d$, where
$r_\ve=1/R_\ve$.

We do a next transformation. Let $z=\ve t$. Then we have $\f(\ve
t,\ve)=\ve g(t,\ve)$ where
\[
\lb{asr3}
\begin{aligned}
&  g(t,\ve)=-A-\ve B-(\a-\ve \b-\ve^2 \b_1)t+\ve t^2\f_1(\ve t,\ve)
=-g_o(t)+\ve g_1(\ve t,\ve),
\\
& g_o(t)=A+t\a, \qq g_1(\ve t,\ve)=-B+t\b+\ve t\b_1+t^2\f_1(\ve
t,\ve).
\end{aligned}
\]
Consider the zeros of the determinant $f_o(z)=\det g_o(z)$. We have
$$
f_o(t)=\det\cT_{p-1}[2a_p'+tb_p]\gt_{p}= \det \cT_{p}\det
[2a_p'+tb_p],
$$
where $\det [2a_p'+tb_p], t\in \C$ is the real polynomial of the
order $d$. This polynomial has only real zeros since the matrix
$a_p', b_p$ are self-adjoint. Let $t_j^o, j\in \N_d$ be zeros of
$f_o$ with the multiplicity $d_j^o$  for some $m$ and $\sum_1^m
d_j^o=d$. Then the determinant $\det g(t,\ve)$ has the form
$$
f(t,\ve)=\det g(t,\ve)=(-1)^df_o(t)+\ve f_1(t,\ve),
$$
where $f_1(t,\ve)$ is an d-th degree polynomial in $t$ and whose
coefficients are all polynomial in of $\ve$. We apply Theorem
\ref{Tp1} to the function $f(t,\ve)$. Due to Theorem \ref{Tp1} for
any $\ve >0$ and each $j\in \N_m$ there exist  $d_j^o$ zeros
$t_{j,s}(\ve)$ of $f(t,\ve)$ in each disc $\dD_{r_\ve}(t_j^o)$ for
some $r_\ve>0$ small  enough. The zeros $t_{j,s}(\ve)$ are
multivalued analytic functions in $\ve$ (not necessary distinct)
with convergent Puiseux series \er{Ps} (Taylor series in
$\ve^{1\/d}$). Thus we have $t_{j,s}(\ve)=t_j^o+O(\ve^{1\/d_j})$.
Then using the transformations $k={1\/z}$ and $z=\ve t$, we obtain
$$
k_{j,s}(\ve)=1/(\ve t_{j,s}(\ve))= (\ve)^{-1}
(t_j^o+O(\ve^{1\/d_j}))^{-1}=(\ve)^{-1}(\t_j^o+O(\ve^{1\/d_j})),
\qqq \where \qq \t_j^o=1/t_j^o.
$$
 \BBox

\subsection{1-dim case}
We consider Jacobi operators $H$ acting on $y=(y_x)_1^\iy\in
\ell^2(\N)$ and given by
\[
\lb{jd1} (Hy)_x=a_{x-1}y_{x-1}+a_xy_{x+1}+b_xy_x, \qq x\in \N,
\]
where $a_x>0, b_x\in \R$ and formally $y_0=0$. In the case $d=1$ the
set $\J_q$ consists of finitely supported sequences $a_x-1, b_x$,
$x\in\N$ and the condition \er{qs} has the form
\[
 \lb{qs1}
 a_{x-p}>0,\
a_{x}=1,\qq b_{x} = 0,\qq \forall \ x>p \qq
 \ca a_{p}=1, b_{p}\ne0\ \ & {\rm if}\ q= 2p-1
 \\
    a_{p}\ne 1,  \ \ & {\rm if}  \ q=2p   \ac .
\]
In order to describe eigenvalues and resonances  of $H$ for this
case we define the set of all possible locations of the eigenvalues
and resonances from \cite{KL23}.

\no {\bf  Definition R}. {\it Let $\cR_q\ss\C^q\sm \{0\}, q\in \N$
be a set of vectors $r=(r_n)_{1}^{q} $ such that:

\no 1) $ 0<|r_1|\le |r_2|\le.....\le |r_q|$ and if $|r_j|=|r_{j+1}|$
for some $j$, then $0\le \arg r_j\le \arg r_{j+1}<2\pi $.  A
polynomial $F(k)=\prod_1^{q}(k-r_n) $ is real on the real line.

\no 2) All zeros of $F$ in $\ol{\dD}$ are real and simple. Denote
the zeros of $F$ in $\dD$ as $k_j,\; j \in \gn_{\bu}$ for some
$\gn_-, \gn_+\ge 0$ arranged by
$$
-1<k_{-\gn_-}<..<k_{-1}<0<k_{1}<...<k_{\gn_+}<1,
$$
and they satisfy
\[
\lb{R1}
 F(1/{k_j})\ne 0 \ \ \  \forall \ j \in \gn_{\bu}=
\{-\gn_-,\dots, -1, 1, \dots, \gn_+\}.
\]
\no 3) $F$ has an odd number $\ge 1$ of zeros on each  interval
$(k_{j}^{-1}, k_{j-1}^{-1}),\; j = \gn_{\bu}\setminus\{-1\}$ and an
even number $\ge 0$ of zeros on each of intervals $(1,
k_{\gn_+}^{-1})$ and $(k_{-\gn_-}^{-1},-1)$.

\no 4) If $q=2p\ge 2$, then ${1\/F(0)}<1$ or it is equivalent real
$F(0)\notin [0,1]$. }

Here the condition ${1\/f(0)}<1$ in 4) is equivalent to $a_p\ne 1$
at $q=2p$.

 We sometimes write $\p_n(V), r_{n} (V),.$. instead of
$\psi_n,r_{n}, .$. when several potentials are being dealt with. Now
we define the mapping $\gr: \J_q\to \cR_q,\; q\in\N$ by
\[
\lb{rmap}
 V \to \gr(V)=(r_n(V))_1^q\in \cR_q,
 \]
where $(r_n)_1^q$ is the sequence of zeros of the corresponding Jost
function $\p(k,V)$. Recall results about inverse resonance
scattering for Jacobi operators (see e.g., \cite{KL23}) with $d=1$.

\begin{theorem}
\lb{T31} i)    Each mapping $\gr :  \J_q\to \cR_q,\; q\in \N$ given
by \er{rmap}  is a bijection between $\J_q$ and $\cR_q$. There is an
algorithm to recover the perturbation $V$ from $\gr(V)$.

\no ii) Let $S_j$ be the S-matrix for two Jacobi operators $H_j$
with perturbations $V_j\in \J_q,\; q\geq 1$ for $j=1,2$. Let
$(\z_j)_1^{q}$ be a sequence of distinct points on $\C_-\bigcap
\{|z| = 1\}$ such that $S_1(\z_j) =S_2(\z_j)$ for all $j\in\N_{q}$.
Then we have $V_1 =V_2$.
\end{theorem}

\no {\bf Remark.} This theorem solves the inverse resonance problem
for Jacobi operators on the half-lattice with finitely-supported
perturbations.

\medskip

From Theorem \ref{T31} we obtain the following results.

\begin{corollary}
\lb{Td1}

 i)  Let finite sequence $r_j\in (-1,0)\cup (0,1), j\in N_m$ for
some $m\ge 0$. Then there exists a finite perturbation $V\in \J_q$
for some $q\ge 1$ such that each  $r_j, j\in N_m$ is an eigenvalue
for the operator $H$.

\no ii) Let a vector $r=(r_n)_{1}^{q}\in \cR_q$ for some $q>1$. Then
there exists a finite perturbation $V\in \J_q$ for some $q\ge 1$
such that each  $r_j\in (-1,1)$ is an eigenvalue for the operator
$H$ and other points $r_j$ are resonances for the operator $H$.

\no iii) The multiplicity of a resonance can be any number for some
specific perturbation $V\in \J_q$.
 \end{corollary}

\no {\bf Example, separation of variables.} Consider Jacobi
operators $H$ on $\mH=\ell^2(\N,\C^d)$ given by \er{Jo}. Let all
matrices $a_x, b_x$ are diagonal and have the from:
$$
a_x=\diag\{a_{1,x},a_{2,x},...,a_{d,x}\},\qqq
b_x=\diag\{b_{1,x},b_{2,x},...,b_{d,x}\}\qqq \forall x\in \N,
$$
for some $(a_{\t,x},b_{\t,x})\in \R_+\ts \R, (\t,x)\in\N_d\ts \N$.
Then the operator $H$ on $\mH=\ell^2(\N,\C^d)$ is an orthogonal sum
of a scalar Jacobi operator $H_\t$ acting on $\ell^2(\N)$ and given
by
\[
\begin{aligned}
\label{u1x}  H=\os_1^d H_\t, \qqq  (H_\t
y)_x=a_{\t,x}y_{x+1}+a_{\t,x-1}y_{x-1}+b_{\t,x}y_x,\qq x\in\N, \qq
y_0=0,
\end{aligned}
\]
where $y=(y_x)_{x\in \N}\in \ell^2(\N)$. In this case each
$f_{\t,x}(k)$ is the standard Jost solution for $H_j $, and
$\p_{\t}(k)=f_{\t,0}(k)$ is the corresponding Jost function. The
properties of $\p_\t$ is well known, see e.g, \cite{KL23} and
\cite{BNW05}, \cite{DZ19}, \cite{IK10}--\cite{IK12}, \cite{K11},
\cite{KL23}, \cite{MNSW12}. In the scalar case the condition \er{qs}
has the form  \er{qs1}. We assume that each sequence $(a_{\t,x},
b_{\t,x}, x\in\N), \t\in \N_d$ satisfies condition \er{qs1}. Note
that the Jost function for the operator $H$ has the for
$\p(k)=\prod_1^d\p_\t(k)$.

\no {\bf Proof of Theorem \ref{T3}.} We consider Jacobi operators
$H$ on $\ell^2(\N)$  given by \er{u1x}, i.e., $ H=\os_1^d H_\t$,
where $(H_\t y)_x=a_{\t,x}y_{x+1}+a_{\t,x-1}y_{x-1}+b_{\t,x}y_x$,
$x\in\N$, and the numbers $(a_{\t,x},b_{\t,x})\in \R_+\ts \R,
(\t,x)\in\N_d\ts \N$ satisfy \er{qs1}. Thus the corresponding
operator $V\in \J_q$.

\no i)  Let  a finite sequence $k_j\in (-1,0)\cup (0,1), j\in N_m$
for some $m\ge 0$, where each $r_j$ has the multiplicity $d_j\le d$.
We take the Jacobi operator $H$ acting on $\ell^2(\N,\C^d)$ given by
\er{u1x}. Due to Theorem \ref{T31} we can assume that $H_1$ has
simple eigenvalues $k_j\in (-1,0)\cup (0,1), j\in N_m$ for some
operator $V_1$ acting on $\ell^2(\N)$ and given by
\[
\lb{V1} (V_1y)_x=(a_{1,x-1}-1)y_{x-1}+b_{1,x}y_x+ (a_{1,x}-1)y_{x+1}
\]
from $\J_q$ for any $q\ge 2m-2$. After this remaining numbers $k_j$
have multiplicity $d_j-1$. From this remaining sequence  we again
due to Theorem \ref{T31} we assume that $H_2$ has these numbers as
simple eigenvalues for some operator $V_2$ acting on $\ell^2(\N)$
and having the form \er{V1} from $\J_q$ for any $q\ge 2m-2$.
Repeating this procedure we  construct the operators $H_1,...,H_d$
and a finite perturbation $V=\os_1^d H_\t\in \J_q$ for some $q\ge
2m-2$ such that each $k_j, j\in N_m$ is an eigenvalue for the
operator $H$.

\no ii) Let a vector $r=(r_n)_{1}^{m}\in \cR_m$ for some $m>0$. Then
there exists a finite perturbation $V\in \J_q$ for some $q\ge 1$
such that each  $r_j\in (-1,1)$ is an eigenvalue for the operator
$H$ and other points $r_j$ are resonances for the operator $H$.

 Let we have a sequence  $\mR=\{r_n\in \C\sm \dD, n\in\N_{qd}\}$ for
some $q\ge 1$, where a number of the points $r_n=1$ (or $r_n=-1$)
$\le d$. From the main sequence  $\mR$ we take the first sequence
$\mR_1=\{\vr_n, n\in\N_{q}\}$ such that $\vr_1=-1$ (if there exists
such point in $\mR$) and $\vr_q=1$ (if there exists such point in
$\mR$) and $\vr_n, n=2,...,q-1$ are any points from $\mR$. For this
sequence $\mR_1$ due to Theorem \ref{T31} there exists a scalar
Jacobi operator $H_1$ with a perturbation $V_1\in \J_q$ having the
form \er{V1} such that each point from this small sequence $\mR_1$
is a resonance for the scalar Jacobi operator $H_1$. From the
remaining sequence $\mR\sm \mR_1$  due to Theorem \ref{T31} we again
take a new sequence  $\mR_2=\{\vr_n, n\in\N_{q}, \vr_n\in \mR\sm
\mR_1\}$ such that there exists a scalar Jacobi operator $H_2$ with
a perturbation $V_2\in \J_q$ having the form \er{V1}. Here each
point from this second small sequence $\mR_2$ is a resonance for the
scalar Jacobi operator $H_2$. We repeat this procedure and we obtain
$H=\os_1^d H_\t$, where each $H_\t$ is the scalar Jacobi operator
with a perturbation $V_\t\in \J_q$ and each point of the sequence
$\mR=\{r_n, n\in\N_{qd}\}$ is a  resonance for some operator $H_\t$.
Here the perturbation $V=\os_1^d V_\t$ belongs to $\J_q$ and we have
used the identity $\p(k)=\prod_1^d\p_\t(k)$. The statement iii)
follows from ii). \BBox

\no {\bf Proof of Theorem \ref{T4}.} i) Assume that Jacobi operators
$H_s=H_o+V_s, s=1,2$ with different perturbations $V_s\in \J_q$ has
the same S-matrix $S(k)$ for all $ k\in \S$. Then the identity
\er{sm2} gives that $V_1, V_2$ have  the same Weyl-Titchmarsh
matrix-valued function $M_{p+1}$:
$$
M_{p+1}(\l)={S(k)-\1 \ol k^{2p+2}\/kS(k)-\ol k^{(2p+3)}} \qqq k\in
\S,
$$
which has a meromorphic extension on the whole plane. Then Theorem
\ref{Thm} gives that $V_1=V_2$.

\no ii) Let $\p=\p(k)$ be the Jost  matrix for the Jacobi operator
$H_s=H_o+V$ with the  perturbation $V\in \J_q$.  Let $\z_j\in \C,\
j\in\N_{q+1}$ be a sequence of distinct numbers. Then due to
\er{Ja1} the matrix-valued polynomial $\p=\p(k)$ has asymptotics
$$
\p(k)=C k^{q}+O(k^{q-1})\qq \as \qqq k\to \iy,\qq \det C\ne0,
$$
for some matrix $C$.  Thus each component of $\p(k)$ is a polynomial
of degree $\le q$ and is given at of distinct numbers $\z_j\in \C,\
j\in\N_{q+1}$. Then we recover it uniquely and then the Jost  matrix
$\p(k)$ is also recovered uniquely. Thus the statement i) and the
identity $S(k)=\p(k)^{-1}\p(1/k)$ yields the Jost function $\p(k)$
determine the perturbation $V$ uniquely. \BBox

We determine the derivative of the phase shift $\x$.

\begin{lemma}
\lb{lpr}  Let  $(r_j)_1^{q}$  be zeros of the Jost function $\gf(k)$
for some perturbation $V \in \J_{q},$ where $q\in \N$. Then for
$k\in \S$ its phase shift function $\x(k)$ satisfies
 \[
 \lb{ss}
\xi'(k)={1\/2ik}
 \sum_{j = 1}^q \frac{2- \Re (\ol k r_j)}{|k- r_j|^2}, \qqq k\in \S.
\]
\end{lemma}
\no {\bf Proof.} Taking the derivative of the both parts of the
identity for the S-matrix
\[\lb{s}
e^{-2i\xi(k)}=\det S(k) =\gf(k)^{-1}\gf(\Bk), \qqq |k| = 1, \qq \Im
k<0,
\]
where $\gf(k)=  \det \p(k)$ and $\Bk=k^{-1}$,  we obtain:
\[
-2ie^{-2i\xi(k)}\xi'(k)=-\frac{\gf'(\Bk)}{k^2\gf(k)}
-\frac{\gf(\Bk)\gf'(k)}{\gf^2(k)}=-\det S(k)
\rt(\frac{\gf'(\Bk)}{k^2\gf(\Bk)}  + \frac{\gf'(k)}{\gf(k)}  \rt),
\]
and
\[
2i\xi'(k)=\frac{\gf'(\Bk)}{k^2\gf(\Bk)}  + \frac{\gf'(k)}{\gf(k)}.
\]
Here we have used that $\psi(k)$ is a polynomial with real
coefficients, for every complex  root $r_j$ there is also a root
$\bar{r_j}$. From Theorem \ref{T1} we have $\p(k) = C\Pi_{j = 1}^q
(k- r_j), C\ne 0$, which implies
$$
\frac{\gf'(k)}{\gf(k)} = \sum_{j = 1}^q\frac{1}{k- r_j}
 =\sum_{j = 1}^q\frac{\bar{k} - \bar{r_j}}{|k - r_j|^2},
 \qqq \frac{\gf'(\Bk)}{\gf(\Bk)} = \sum_{j = 1}^q\frac{1}{{1\/k}- r_j}
 =\sum_{j = 1}^q\frac{1}{\ol k- \ol r_j}
$$
and
$$
\begin{aligned}
2i\xi'(k)=\frac{\gf'(1/k)}{k^2\gf(1/k)}  + \frac{\gf'(k)}{\gf(k)} =
\sum_{j = 1}^q\frac{\bar{k} - \bar{r_j}+
 \bar{k}^2(k- r_j)}{|k- r_j|^2} =
 \sum_{j = 1}^q \frac{2- \Re (\ol k r_j)}{k|k- r_j|^2}.
\end{aligned}
$$
  \BBox

\section {\lb{Sec5}Fredholm determinants}
\setcounter{equation}{0}

\subsection {Laplacians on $\Z$} Consider the Laplacian $\D$ on the lattice $\Z$
given by $(\D y)_x=y_{x+1}+y_{x-1}, x\in \Z$.
 Let us introduce the shift operator on the lattice $\Z$ (and $\N$)
 by
\begin{equation}
(\pa y)_x=y_{x+1}, \quad (\pa^* y)_x=y_{x-1}.
\end{equation}
 The spectral properties of $\D$ is easier to  describe by passing
to its unitary transformation by the Fourier series. Let $ {\T}=
{\R}/(2\pi \Z)= [-\pi,\pi]$ be the flat torus and ${\cU}$ be the
unitary operator from $\ell^{2}({\Z})$ to $L^{2}({\T}),$ given  by
\begin{equation}
\tes \wh y(t)=({\cU} y)(t) =  \sum_{x\in{\Z}}y_x {e^{-it\cdot
x}\/(2\pi)^{1\/2}},\qq t\in \T. \nonumber
\end{equation}
The shift operator and the Laplacian on $L^{2}({\T})$ are rewritten
as $\wh \pa := \mathcal U\,\pa\,\mathcal U^{\ast} = e^{it}$ and
\begin{equation}
\label{Ls} \wh \D={\mathcal U}\, \D\, {\mathcal U}^{\ast}={2}\cos
t,\qqq t\in \T.
\end{equation}
Let $\gR(k)=(\D-\l(k))^{-1}$ and let $\gR(x-x',k), (x,x',k)\in
\Z^2\ts \dD$ be its kernel.
 Then the kernel $\gR(x,k)$ has the following representation
\begin{equation}
\lb{r1} \gR(x, k) = \frac{k^{|x|}}{k - 1/k}
=-\frac{k^{|x|}}{\sqrt{\lambda^2-4}},\quad {\rm for} \quad (x,k)\in
{\Z}\times {\dD},
\end{equation}
and $\gR(x,k)$ has a meromorphic continuation from ${\dD}$ into
${\C}$, see e.g. \cite{IK12}. The proof is very simple. Indeed,
residue calculus imply
$$
\gR(x, k)=\frac{1}{2\pi}\int_0^{2\pi}\frac{e^{-ixt}dt}{2\cos t-k -
{1\/k}} =\frac{1}{2\pi i}\int_{|w|=1}\frac{w^{-x}dw}{(w-k)(w-
{1\/k})} = \frac{k^{|x|}}{k - {1\/k}}.
$$
By \er{kL}, we have ${1\/k} = {\lambda+\sqrt{\lambda^2-4}\/2}$ and
$k -{1\/k}= - \sqrt{\lambda^2-4}$, which proves \er{r1}.

\begin{lemma} \lb{Tr1}
i)  The kernel  $R_o(x,x',k), (x,x',k)\in {\N^2}\times {\dD}$ of the
resolvent $R_o(k)=(H_o-\l(k))^{-1}$ on $\ell^2(\N)$ is a polynomial
and  has the form
\begin{equation}
\lb{r2}
\begin{aligned}
 R_o(x,x',k)= \frac{k^{|x-x'|}-k^{x+x'}}{k - 1/k}=k^{|x-x'|}
\frac{1-k^{2w}}{k - 1/k},
 \end{aligned}
\end{equation}
\begin{equation}
\lb{r2z}
\begin{aligned}
R_o(x,x',k)=-k^{1+|x-x'|}(1+k^2+k^4+...+k^{2(w-1)}),\qq
w=\ca x', \ x\ge x'\\
x, \ x'>x\ac,
 \end{aligned}
\end{equation}
and it  has an analytic continuation from ${\dD}$ into ${\C}$.

\no ii) The following estimates hold true for $x\le x'$:
\[
\lb{r3} |R_o(x,x',k)|\le C|k|^{|x-x'|+1},\qq C=x\qq  or \qq
C={2\/|k^2-1|},\qq  |k|\le 1,
\]
\[
\lb{r4} |R_o(x,x',k)|\le C_2|k|^{x'+x-1}, \qq C_2=\ca  x  ,
\qqq \ if \qq |k|\le 2\\
                  4   ,\qqq \ if \qq |k|\ge 2                 \ac .
\]
\end{lemma}

\no {\bf Proof.} i) Let $L^2(\T)=L^2_{ev}(\T)\os L^2_{od}(\T)$,
where $L^2_{od}(\T)$ is the set of odd functions $f(-t)=-f(t)$ and
$L^2_{ev}(\T)$ is the set of even functions $f(-t)=f(t)$ for all
$t\in [-\pi, \pi]$. Introduce the similar decomposition
$\ell^2(\Z)=\ell^2_{ev}(\Z)\os \ell^2_{od}(\Z)$, where is
$\ell^2_{ev}(\Z)$ the set of odd functions $y_{-x}=-y_{x}$ and
$\ell^2_{od}(\Z)$ is the set of even functions $y_{-x}=y_{x}$ for
all $x\in \Z$. From \er{Ls} we obtain
\[
\lb{DH}
\begin{aligned}
 \D \ell^2_{s}(\Z)\ss \ell^2_{s}(\Z), \qq s\in \{od, ev\},
 \\
 \D y=H_oy\qqq \forall \qq y\in \ell^2(\N),\qq y_0=0.
 \end{aligned}
\]
For each $h\in \ell^2(\N)$ we define an odd function $h^{od}\in
\ell^2_{od}(\Z)$ by
$$
 h_x^{od}=\ca  h_x\qq if \qq x\in \N\\
              -h_x\qq if \qq -x\in \N \ac\qq  \and \qq \wt h_0=0.
$$
Let $\a={1\/k - 1/k}$. We have
$$
g_x=\a \sum_{j\in \N}(k^{|x-j|}-k^{|x+j|})h_j=\a \sum_{j\in
\Z}k^{|x-j|}h_j^{od}=(\gR(\l(k))h^{od})_x=g_x^{od},\qq \forall \
x\in \N,
$$
and $(\D-\l)g^{od}=h^{od}$, which yields
$((\D-\l)g^{od})_x=((H_o-\l) g)_x=h_x^{od}=h_x$ for $x\in\N$, since
$h_0^{od}=g_o^{od}=0$. From here we obtain \er{r2}, which yields
\er{r2z}.

\no ii) Using \er{r2} we obtain
\begin{equation}
\lb{r5} R_o(x,x',k)= \frac{k^{|x-x'|+1}}{k^2 - 1} (1-k^{2w}),\qq
2w=x+x'-|x-x'|=\ca 2x' \qq if \qq x\ge x'\\
               2x \qq if \qq x\le x' \ac,
\end{equation}
which yields \er{r2}.   Let $x\le x'$ and $\ve=k^2$. Then we get
\begin{equation}
\lb{r6} |1-k^{2x}|\le|1-\ve ||1+\ve+....+\ve^{x-1}|\le |1-k^2|x\qqq
\forall\qq k\in \C.
\end{equation}
Let $|k|\le 1$. The estimate for the case $C={2\/|k^2-1|}$ follows
from  \er{r2}. For the case $C=x$,  \er{r6} implies $|1-k^{2x}|\le
|1-k^2|x$ and thus \er{r5} yields
\[
\begin{aligned}
|R_o(x,x',k)|\le \frac{|k|^{|x-x'|+1}}{|k^2 - 1|} |1-k^{2w}|\le
x|k|^{|x-x'|+1}.
\end{aligned}
\]
If $|k|\le 2$, then   \er{r6} implies $|1-k^{2x}|\le |1-k^2|x
k^{2(x-1)}$ and thus \er{r5} yields
\[
\begin{aligned}
 |R_o(x,x',k)|\le \frac{|k|^{|x'-x|+1}}{|k^2 - 1|} |1-k^\b|\le
x|k|^{x'-x+1+2(x-1)}=x|k|^{x'+x-1}.
\end{aligned}
\]
If $|k|\ge 2$, then \er{r2} yields
$$
|R_o(x,x',k)|\le \frac{2|k|^{x+x'+1}}{|k^2 - 1|}\le 4|k|^{x+x'-1}.
$$
Collecting the last two estimates we obtain \er{r4}. \BBox

We determine the kernel $R_{x,x'}(k), x,x'\in \N$ of the resolvent
$(H-\l)^{-1}, \l=k+{1\/k}, k\in \dD$, see e.,g., p. 152 from
\cite{CGR05}.

\begin{proposition}
\lb{Trr1} Let $V\in \J_q$ and  $\l=k+{1\/k}, k\in \dD$. Then

\no i)  The resolvent $R(k)=(H-\l)^{-1}, \l\in \L_1$ on $g\in
\ell^2(\N,\C^d)$ is given by
\[
\lb{dR}
\begin{aligned}
y=R(k)g=\gF(k)u(k)+\vp(k)v(k),
\\
u_x(k)=\sum_{t=1}^x \wt\vp_t(k)g_t,\qq v_x(k)=\sum_{t=x+1}^\iy \wt
\gF_t(k) g_t,
\end{aligned}
\]
ii) Let $R_{x,t}(k), x,t\in \N$ be the kernel of the resolvent
$R(k), k\in \dD$. For each $x,t\in \N$ the matrix-valued function
$R_{x,t}(k)$  has a meromorphic extension from the unit disc in the
whole complex  plane with poles at zeros of the Jost function
$\gf(k)=\det \p(k)$.

\end{proposition}

\no {\bf Proof.} i) Consider the equation $(H-\l)y=g, \l\in \L$. We
construct the solution
$$
 y_x=\gF_xu_x+\vp_xv_x,\qq u_0=0,
 $$
in terms of some unknown matrix-functions $u_x, v_x$. Assume that
components $\z_x=u_{x+1}-u_x, \ \g_x=v_{x+1}-v_x$ satisfy
\[
\lb{ac1}
\begin{aligned}
\gF_{x}\z_{x}+\vp_{x} \g_{x}=0,
\\
a_{x}\Big[\gF_{x+1}\z_{x}+\vp_{x+1} \g_{x}\Big]=g_x,
\end{aligned}
\]
for $x\in \N$.  We compute  $(Hy)_x =
 a_{x-1}y_{x-1}+ a_{x}y_{x+1}+b_x y_k$
for $ y_x=\gF_xu_x+\vp_xv_x. $
 The components $y_{x-1}, y_{x+1}$ in terms of $\z_x, \g_x$ have the forms
\[
\lb{rt}
\begin{aligned}
 & y_{x-1}= \gF_{x-1}u_{x-1}+\vp_{x-1}
v_{x-1}=(\gF_{x-1}u_{x}+\vp_{x-1}
v_{x})+(\gF_{x-1}\z_{x-1}+\vp_{x-1} \g_{x-1}),
\\
& y_{x+1}=\gF_{x+1}u_{x+1}+\vp_{x+1}v_{x+1}
 =(\gF_{x+1}u_{x}+\vp_{x+1}
v_{x})+(\gF_{x+1}\z_{x}+\vp_{x+1} \g_{x}).
\end{aligned}
\]
 Substituting \er{rt} into $(Hy)_x$ we compute
$$
\begin{aligned}
(Hy)_x=a_{x-1}\Big[\gF_{x-1}u_{x}+\vp_{x-1} v_{x} \Big] + a_{x}\Big[
\gF_{x+1}u_{x}+\vp_{x+1} v_{x}\Big]+g_x+b_x y_k=\l y_x+g_x.
\end{aligned}
$$
We need to show that the system \er{ac1} has a unique solution. From
\er{ac1} we obtain
$$
\begin{aligned}
\g_{x}=-\vp_{x}^{-1}\gF_{x}\z_{x} \ \Rightarrow \ G_x\z_x=g_x,\qqq
G_x=a_{x}\gF_{x+1}-a_{x}\vp_{x+1}\vp_{x}^{-1}\gF_{x}.
\end{aligned}
$$
Then from this solution and  \er{wr1}, \er{wt1}  we get
$$
%a_{x}\gF_{x+1}-\vp_{x+1}\vp_{x}^{-1}\gF_{x}
G_x=
a_{x}\gF_{x+1}-\wt\vp_{x}^{-1}\wt\vp_{x+1}a_{x}\gF_{x}=\wt\vp_{x}^{-1}
\Big[\wt\vp_{x}a_{x}\gF_{x+1}-\wt\vp_{x+1}a_{x}\gF_{x}    \Big]=
\wt\vp_{x}^{-1}\{\vp, \gF \}=-\wt\vp_{x}^{-1},
$$
which yields $\z_x=-\wt\vp_{x}g_x$. We need to determine $\g_x$.
From \er{ac1} we obtain
$$
\begin{aligned}
\z_{x}=-\gF_{x}^{-1}\vp_{x}\g_{x} \ \Rightarrow \
a_{x}\Big[\vp_{x+1}-\gF_{x+1}\gF_{x}^{-1}\vp_{x}\Big]\g_x=g_x.
\end{aligned}
$$
Then from \er{ac1}, \er{wt1},   we get
$$
\begin{aligned}
a_{x}\vp_{x+1}-a_{x}\gF_{x+1}\gF_{x}^{-1}\vp_{x}=
 a_{x}\vp_{x+1}-\wt \gF_{x}^{-1}\wt \gF_{x+1}a_{x}\vp_{x}=
\wt \gF_{x}^{-1} \Big[\wt \gF_{x}a_{x}\vp_{x+1}-\wt
\gF_{x+1}a_{x}\vp_{x} \Big]=\wt \gF_{x}^{-1},
\end{aligned}
$$
which yields $\g_x=\wt\gF_x g_x$. We compute $u_x, v_x$. We have
$$
u_{x}=\z_{x-1}+u_{x-1}=\z_{x-1}+\z_{x-2}+u_{x-3}
=....=\z_{x-1}+\z_{x-2}+...+\z_{1}.
$$
Let $r_x(n)=-\sum_x^{n-1}\g_j$.  Similar arguments give
$$
v_x=-\g_x+v_{x+1}=-\g_x-\g_{x+1}+v_{x+2}=-\g_x-\g_{x+1}-\g_{x+2}+v_{x+2}
=
r_x(n)+v_{n+2}.
$$
The sum $\sum_x^{\iy}\g_j$ is absolutely convergence since
$\g_j=\gF_j g_j$. These identities and $\z_x=-\wt\vp_{x}g_x$ and
$\g_x=\wt\gF_x g_x$ imply \er{dR}.

\no ii) From i) and the properties of the functions $\vp, F=\vt+\vp
m_+, m_+=f_1f_0^{-1}$ we deduce that for each $x,t\in \N$ the
matrix-valued function $R_{x,t}(k)$  has a meromorphic extension
from the unit disc in the plane $\C$ with poles at zeros of the Jost
function $\gf(k)=\det \p(k)$. \BBox

\begin{lemma}
\lb{Tr2}  Let $\|\cdot\|_{\B_2}$ be the Hilbert-Schmidt norm on
$\ell^2({\N})$ and $k\in \dD$. Then we have
\begin{equation}
\lb{rr1} \|f R_o(k)g\|_{\B_2}\le \frac{|k|}{|k^2-1|}\|f\|\|g\|,\qqq
f, g\in \ell^2({\N},\C^d),
\end{equation}
\begin{equation}
\lb{rr2} \|f R_o(k)g\|_{\B_2}\le \sqrt 2 |k|\|f\|_1\|g\|_1, \qqq f,
g\in \ell_1^2({\N},\C^d),
\end{equation}
where  $\|f\|^2=\sum_{x\in \N} |f_{x}|^2 $ and
                       $\|f\|_1^2=\sum_{x\in \N}x^2 |f_{x}|^2 $.
\end{lemma}

\no {\bf Proof.} From \er{r3} we obtain
$$
\|f R_o(k)g\|_{\B_2}^2= \sum_{x,x'\in \N}
|f_{x}|^2{|k^{|x-x'|}-k^{|x+x'|}|^2\/|k - 1/k|^2} |g_{x'}|^2 \le
{1\/|k - 1/k|^2} \sum_{x,x'\in \N} |f_{x}|^2 |g_{x'}|^2,
$$
which yields \er{rr1}.  From \er{r3} we obtain
$$
\begin{aligned}
\|f R_o(k)g\|_{\B_2}^2= \sum_{x,x'\in \N}
|f_{x}|^2{|k^{|x-x'|}-k^{|x+x'|^2}|\/|k - 1/k|^2}^2 |g_{x'}|^2
\\
\le |k|^2 \sum_{x>x'} {x'}^2 |f_{x}|^2 |g_{x'}|^2+ |k|^2 \sum_{x<x'}
x^2 |f_{x}|^2 |g_{x'}|^2\le 2|k|^2 \|f\|_1^2\|g\|_1^2,
\end{aligned}
$$
which yields \er{rr2}. \BBox

Recall that we have defined $(\pa f)_x=f_{x-1},\qq (\pa^*
f)_x=f_{x+1}, x\ge 1$ acting on $\ell^2(\N)$ by
$$
\pa f =(0,f_1,f_2,f_3,...),\qqq \pa^* f =(f_2,f_3,f_4,...).
$$
We rewrite $H_c$ and $V_c$ in terms of $\pa$
$$
H_c=\pa a+b+s\pa^*=H_o+V_c, \qqq  V=\pa (a-\1)+b+(s-\1)\pa^*.
$$
We need simple properties of the operator $\pa$:
\[
\lb{trp1}
\begin{aligned}
 &  \pa^*\pa =\1,\qqq \pa\pa^* f=(0,f_2,f_3,f_4,...),
\\
& \Tr {\pa}^n A =0,\qqq \forall \qq (n,A)\in \N\ts\ell^1(\N,\M_d)
\end{aligned}
\]
 Define Fredholm determinants
$\cD_c(\l),  D_c(k)$ and $\cD(\l),  D(k),$ where $\l=k+{1\/k}\in
\L$ by
\[
\lb{DD1}
\begin{aligned}
\cD_c(\l)=\det (I+Y_c(\l)), \qq D_c(k)=\cD_c(\l(k)),\qq
Y_c(\l)=V_cR_o(\l)\c_c,
\\
\cD(\l)=\det (I+Y(\l)), \qq D(k)=\cD(\l(k)),\qq Y(\l)=VR_o(\l)\c,
\end{aligned}
\]
where $\c_c$ (and $\c$) is the characteristic function of the set
$\supp (|a_x-\1|+|s_x-\1|+|b_x|)$ (and $\supp (|a_x-\1|+|b_x|)$),
where $|\cdot|$ is a matrix norm in $\C^d$. The determinants are
well defined since $V_c, V$ are finitely supported perturbations and
thus $Y_c(\l), Y(\l)\in \cB_1$.

\begin{lemma}
\lb{Tr3}
i)   Let functions  $F, G\in \ell^2(\N, \M_d)$ be finitely
supported. Then  operator-valued functions
$\cT(k)=(k^2-1)FR_o(k)G\in \cB_1$ for any $k\in \dD$, the function
$\cT:\dD\to \cB_1$ is analytic in $\dD$ and is a polynomial in $\C$
with finite dimensional coefficients.

\no ii)  Let   $V\in \J_q^c$. Then the operator
$\cT(k)=(k^2-1)Y_c(k)\in \B_1$ for any $k\in \dD$ and the
operator-valued  function $\cT:\dD\to \B_1$ is analytic in $\dD$ and
has an analytic extension onto the whole plane, which is a
polynomial in $\C$ with finite dimensional coefficients.

\no iii) The Fredholm determinant $D_c(k), k\in \dD$ is analytic in
$\dD$  and has a meromorphic continuation from ${\dD}$ into the
whole plane ${\C}$, which is a polynomial in $\C$.

\no iv)  The set $\{k \in \dD: D_c(k) = 0\}$ is finite and coincides
with $\sigma_d(H_c)$.

\end{lemma}

\no {\bf Proof.} i) Let $\a=|x-x'|$. From \er{r2} we obtain
\begin{equation}
\lb{rq3} R_o(x,x',k)=k^{1+2\a} \frac{1-k^{2w}}{k^2 -
1}=k^{1+2\a}(1+k^2+k^4+...+k^{2(w-1)}),
\end{equation}
and each $R_o(x,x',k), x,x'\in \N$ has a meromorphic continuation
from ${\dD}$ into ${\C}$.

\no Results ii) and iii) follows from i). The statement iv) is
well-known. \BBox

We need results from \cite{IK12}  about Fredholm determinants.

\begin{lemma}\lb{TIK12}

Let an operator $H = H_o + V$ act on the Hilbert space $\mathcal K$,
where the operators $H_o=H_o^*, V$ are bounded and the operator $\
V$ is trace class. Let $\cD(\lambda)=\det (I+VR_o(\lambda))$, where
$R_o(\lambda) = (H_o - \lambda)^{-1}, \lambda \in \Lambda := \C\sm
\s(H_o)$. Then

 \no i) $\ D(\lambda)$ is analytic in $\Lambda$. Moreover
 \begin{equation}
 \label{prD}
\ol\cD(\lambda)=\cD(\ol\lambda),\qqq \l\in \L,\qq {\rm if}\qq V^*=V,
\end{equation}
\begin{equation}
\cD(\lambda)=1+O(1/\lambda) \quad \as \quad |\lambda|\to {\infty},
\label{S6Dlambdalarge}
\end{equation}
\begin{equation}
\begin{aligned}
\label{d1} & \log \cD(\lambda)= -
\sum_{n=1}^{\infty}\frac{(-1)^n}{n}{\rm
Tr}\,\left(VR_o(\lambda)\right)^n=-\sum _{n \geq 1}\frac{F_n}{n\
\lambda^n},
\\
&F_1={\rm Tr}\, V, \quad F_2={\rm Tr}\,(2VH_o+V^2),\ F_n={\rm
Tr}\,(H^n-H_o^n), \quad n\ge 1,
\end{aligned}
\end{equation}
where the series are absolutely convergent for $|\lambda|
> r_0$ and the radius $r_0 >0$ is large enough.

\no ii) The set $\{\lambda \in \Lambda\, ; \, \cD(\lambda) = 0\}$ is
discrete and coincides with $\sigma_d(H)$.

\end{lemma}

We discuss properties of the Fredholm determinant $D_c(\l)=\det
(I+V_cR_o(\l))$.

\begin{lemma}
\lb{Ttr1}

Let an operator $H_c = H_o + V_c$ act on $\mH$, where $ V_c\in
\J_q^c$. Then

\no i)  The determinant $\cD_c(\l)$ satisfies i) and ii) from Lemma
\ref{TIK12}, where $\L = \C\sm \s(H_o)$.

\no ii)  The function   $\ D_c(k)=\cD_c(\l(k)), \l=k+{1\/k}$ is a
polynomial in $k\in \C$ and satisfies
\begin{equation}
\label{Dk1} \tes \log D_c(k)=-F_1 k-F_2{k^2\/2}+ O(k^3) \quad \as
\quad k\to 0,
\end{equation}
\begin{equation}
\begin{aligned}
\label{Dk2}  F_1=\Tr \sum b_x,\qqq F_2=\sum_{x=1}^{\infty}\Tr
\big(b_x^2+2 (a_{x}s_{x}-\1)\big).
\end{aligned}
\end{equation}

\no iii) The set $\{\lambda \in \Lambda\, ; \, \cD_c(\lambda) = 0\}$
is finite and coincides with $\sigma_d(H_c)$.

\end{lemma}

\no {\bf Proof}. The statement  i) follows from Lemma \ref{TIK12}.

\no ii) Substituting ${1\/\l}={k\/1+k^2}=k+O(k^3)$ into \er{d1} we
obtain
$$ -\log
\cD(\l)={F_1\/\l}+ {F_2\/2\l^2}+{O(1)\/\l^3}={F_1k}+
{F_2k^2\/2}+{O(k^3)}
$$
We compute $F_1, F_2$. Let $\a=a-\1_d, \s=a-\1_d,$ and
$V_o=\pa\a+\s\pa^*$. From  \er{d1}, \er{trp1} we obtain
$$
F_1=\Tr V_c= \Tr (\pa\a+b+ \s\pa^*)=\Tr b,\qqq F_2={\rm
Tr}\,(2V_cH_0+V_c^2).
$$
We compute $F_2$. From \er{trp1} we obtain the first term
$$
\begin{aligned}
&\Tr 2V_cH_o=2\Tr (\pa\a+\s\pa^*+b)(\pa+\pa^*)=2\Tr \big[(\pa\a
+\s\pa^*+b)\pa+  (\pa\a+\s\pa^*+b)\pa^*  \big]
\\
&=2\Tr \big[  \s\pa^*\pa +\pa\a\pa^* \big]=2\Tr \sum_{1}^\iy\big[
\a_x+\s_x \big].
\end{aligned}
$$
Let $V_o:=\pa\a+\s\pa^*$. From \er{trp1} we obtain the second term
$$
\begin{aligned}
\Tr V_c^2= \Tr [V_o +b]^2=\Tr [b^2+2V_o b+ V_o^2]=\Tr [b^2+ V_o^2],
\\
\Tr V_o^2=\Tr (\pa\a+\s\pa^*)^2=\Tr
[\pa\a\pa\a+\s\pa^*\s\pa^*+\pa\a\s\pa^*+\s\pa^* \pa\a]
\\
= \Tr [\pa\a\s\pa^*+\s\pa^* \pa\a]=\sum_{x=1}^\iy \Tr
[\a_{x}\s_{x}+\s_{x}\a_{x}]=2\sum_{x=1}^\iy \Tr \s_{x}\a_{x}
\end{aligned}
$$
Collecting these identities and using $\Tr
(\a_{x}\s_{x}+\a_{x}+\s_{x}=\a_{x}s_x-\1$ we obtain \er{Dk2}.

The statement  iii) follows from Lemma \ref{TIK12}, ii), since
$D_c(k)$ is a polynomial. \BBox

\subsection{Determinants }

Now we discuss the real Jacobi matrices. We prove  the uniqueness.
Let $\mathcal H^2=\mathcal H^2(\dD)$ denote the Hardy class of
functions $f$ which are analytic in the unit disc $\dD$ equipped
with the norm
$$
\|f\|_{\mathcal
H^2}^2=\sup_{1>r>0}{1\/2\pi}\int_0^{2\pi}|f(re^{i\f})|^2d\f<\infty.
$$

\begin{lemma} \label{Ta1}
Let $f_j, j=1,2$, be such that $\log f_j\in \mathcal H^2$, and
satisfy the following conditions:

\noindent i) $f_j(k)=1+o(1)$ as $k\to 0$,

\noindent ii) $\phi_j(t) := \arg f_j(t)$ is continuous in $t\in
{\S}\setminus \{k_o\}$ for some $k_o\in \S$,

\noindent iii) $e^{-2i\phi_1(t)}=e^{-2i\phi_2(t)},$ for a.e. $t\in
{\S}$.

Then $f_1=f_2$ on $\dD$.
\end{lemma}

\no {\bf Proof}. Recall that if $\phi_1(t)=\phi_2(t), $ for a.e.
$t\in {\R}$. Then it is well known that $f_1(z)=f_2(z)$ on $\dD$
(see p.10 \cite{Ko88}). Hence we have only to show
$\phi_1(t)=\phi_2(t), $ for a.e. $t \in {\S}$. Due to ii) and iii)
we deduce that $\f_1(t)=\f_2(t)+\pi N, $ for a.e. $t \in {\S}$ and
for some $N\in \Z$. Hence we have $\log f_1(k)=\log f_2(k)+i\pi N$
and i) yields that $N=0$. \BBox

We discuss properties of $\gf(k)=\det f_0(k)$ and $D(k)$ with real
coefficients.

\begin{lemma}
\lb{Ta2}

Let $V\in \J_q,\; q$ and let  all  $a_x, b_x$ are real. Then we have
\[
\lb{a21}
 \ol\gf(k)=\gf(\ol k),\qqq \forall \ k\in \C.
\]
 Let in addition, $D(1)\ne 0$ or $D(-1)\ne 0$. Then
\[
\lb{a22} \gf=D\det \cT_p.
\]

\end{lemma}
\no{\bf Proof.} Consider $\gf=\det \p$, where the matrix-valued
polynomial $\p(k)$ has matrix valued coefficients, which are real
matrices. Then it gives \er{a21}. In order to show \er{a22} we
define the function $F(k)=\gf(k)/\det \cT_p$. The functions $F(k),
D(k)$ have the same zeros in $\dD$ and
$$
F(k)=1+O(k),\qqq  D(k)=1+O(k)\qq \as \qq k\to 0.
$$
For the zeros $k_j\in (-1,1), j\in \gn_\bu$ of $\gf$ (and
simultaneously $D(k)$)  we define the Blaschke product $B(k), k\in
\C$ by: if $\gn_\bu=0$, then $B=1$ and if $\gn_\bu\ge 1$, then
\[
\lb{B2}
\begin{aligned}
 B(k)=\prod_{j\in \gn_\bu} {|k_j|\/k_j}{(k_j-k)\/(1-\ol k_j k)}.
\end{aligned}
\]
Define the functions $ F_B=F/B $ and $D_B=D/B$. From the
Birman-Krein identity $ \det S(k)={\ol D(k)/D(k)}$, and from
$S(k)=\p(k)^{-1}\p(\ol k), k\in \S$ and \er{a21} we obtain
$$
{\ol D(k)\/D(k)}={\gf(\ol k)\/\gf(k)}={\ol \gf(k)\/\gf(k)}= {\ol
F(k)\/F(k)}.
$$
This yields $ {\ol D_B\/D_B}={\ol F_B\/F_B} $. Thus  due to Lemma
\ref{Ta1} we deduce that $D_B=F_B$ and then $D=F$. \BBox

Now we discuss Jacobi operators $H_c=H_o+V_c$, where $V\in \J_q^c$.
Consider the Jost function $f_x(k,h)$ and the Jost determinant
$\gf(k,h)=\det f_0(k,h)$ for the operator $H_c$, where
$h=(h_x)_1^\iy, h_x=(a_x,s_x,b_x)$. Define determinants and
functions
\[
\lb{dd1}
\begin{aligned}
f_x^\bu(k,h)=(\det A_x) f_x(k,h),\qq \gf^\bu(k,h)=(\det A_x)\det
f_0(k,h), \qq A_x=a_1a_2....a_x.
\end{aligned}
\]

\begin{theorem} \lb{Ta3}

Let $V_c\in \J_q^c,\; q\in\{2p,2p-1\}$. Then each $f_x^\bu(k,h),
x\ge 0$ is a polynomial in $k, a_x^\bu, a_j^\bu, b_j, s_j, j\ge x-1$
and satisfies
\[
\lb{a21z}
\begin{aligned}
& f_0^\bu(k,h)=A_p^\bu-k\cK_{0}^\bu+...,
\\
& A_p^\bu= a_{1}^\bu a_{2}^\bu\cdots a_p^\bu,\qqq
\cK_{0}^\bu=b_1a_{1}^\bu...a_{p}^\bu+a_{1}^\bu b_2 a_{2}^\bu \cdots
a_{p}^\bu+...+ a_{1}^\bu \cdots a_{p-1}^\bu b_{p}a_{p}^\bu.
\end{aligned}
\]
Moreover, the function $\det f_0^\bu(k,h)$ is a    polynomial in $k,
a_x^\bu, s_x, b_x, x\ge 1$ and satisfies
\[
\lb{dea1}
\begin{aligned}
\gf^\bu(k,h)=1+O(k),\qqq \forall \ k\to \iy.
\end{aligned}
\]

\end{theorem}
\no{\bf Proof.} i) Consider the matrix-valued function
$f_x^\bu(k)=\det(A_x)f_x(k,h),x\ge 0$.
 Recall that if $a_x$ is an invertible matrix, then
$\gt_x=a_x^{-1}={1\/\det a_x} a_x^{\bu}$, where $a_x^\bu$ be the
adjugate matrix of a matrix $a_x$ and $\det a_x^\bu=1$.
 Using the equation
$a_{x}f_{x}=(\l-b_{x+1})f_{x+1}-s_{x+1}f_{x+2}$ and $f_x=k^x,\;
x>p$, we obtain $f_p=k^{p}\gt_p$ and $f_p^\bu=k^{p}a_p^\bu$ and the
following cases:\\
 \no $\bu$ if $x=p-1$, then
$$
\begin{aligned}
 &  \tes a_{p-1}f_{p-1}=(\l-b_{p})f_{p}-s_pf_{p+1}
=((k+{1\/k}-b_{p})k^{p}-s_pa_pk^{p+1}) \gt_p
=k^{p-1}(\1-b_pk+\e_pk^{2})\gt_p,
\\
&f_{p-1}^\bu=k^{p-1}a_{p-1}^\bu(\1-b_pk+\e_pk^{2})a_p^\bu, \qqq
\e_x=\1-s_x a_x;
\end{aligned}
$$
$\bu$ if $x=p-2$, then
$$
\begin{aligned}
\tes a_{p-2}f_{p-2}={(\l-b_{p-1})f_{p-1}-s_{p-1}f_{p}}
=\Big[(k+{1\/k}-b_{p-1})\gt_{p-1}(k^{p-1}-b_pk^{p}
+\e_pk^{p+1})-s_{p-1}k^{p} \Big]\gt_p
\\= \Big[\gt_{p-1}k^{p-2}-k^{p-1}(b_{p-1}\gt_{p-1}+
\gt_{p-1}b_{p})+....\Big]\gt_{p},
\end{aligned}
$$
which yields
$$
\begin{aligned}
 f_{p-2}^\bu=\Big[k^{p-2}a_{p-2}^\bu a_{p-1}^\bu a_{p}^\bu-
 \cK_{p-2}^\bu k^{p-1}
+....\Big], \qq \cK_{p-2}^\bu=a_{p-2}^\bu b_{p-1}a_{p-1}^\bu
a_{p}^\bu+a_{p-2}^\bu a_{p-1}^\bu b_{p}a_{p}^\bu.
\end{aligned}
$$
$\bu$ If $x=p-3$, then we have
$$
\begin{aligned}
 f_{p-3}=\gt_{p-3}\Big\{(\l-b_{p-2})f_{p-2}-s_{p-2}f_{p-1}\Big\}
=k^{p-3}\gt_{p-3}\gt_{p-2}\gt_{p-1}\gt_p-\cK_{p-3}k^{p-2} +....,
\end{aligned}
$$
where
$$
\cK_{p-3}= \gt_{p-3}b_{p-2}\gt_{p-2}\gt_{p-1}\gt_p+
\gt_{p-3}\gt_{p-2}
b_{p-1}\gt_{p-1}\gt_p+\gt_{p-3}\gt_{p-2}\gt_{p-1}b_{p}\gt_p,
$$
which yields
$$
\begin{aligned}
& f_{p-3}^\bu=k^{p-3}a_{p-3}^\bu a_{p-2}^\bu a_{p-1}^\bu a_{p}^\bu-
 \cK_{p-3}^\bu k^{p-2}
+....,
\\
&  \cK_{p-3}^\bu=a_{p-3}^\bu b_{p-2}a_{p-2}^\bu a_{p-1}^\bu
a_{p}^\bu+ a_{p-3}^\bu a_{p-2}^\bu  b_{p-1}a_{p-1}^\bu a_{p}^\bu
+a_{p-3}^\bu a_{p-2}^\bu a_{p-1}^\bu b_{p}a_{p}^\bu.
\end{aligned}
$$
Repeating this procedure we obtain  \er{a21z}. Moreover, we deduce
that  $f_x^\bu(k,h), x\ge 0$ is a polynomial in $k, a_x, a_j^\bu,
b_j, s_j, j\ge x$. Furthermore, we obtain that the function $\det
f_0^\bu(k,h)$ is a    polynomial in $k, a_x^\bu, s_x, b_x, x\ge 1$
and satisfies \er{dea1}.
 \BBox

\no {\bf Proof of Theorem \ref{T1x}.} Consider the function
$D_c(k,h)=\det (\1+V_c R_o(k))$. The operator $Y_c(k,h)=
V_cR_o(k)\c_p$ is finite dimensional matrix and then due to \er{r2}
the function $D_c(k,h)$ is a polynomial of $k, h_x$. The functions
$\gf^\bu(k,h), D_c(k,h)$ are polynomials in $k, h_x$. They coincide
at real  $b_x=b_x^*$ and $s_x=a_x>0$. Then they coincide for all $h$
and we have  \er{JfD}. If $V\in \J_q$, then \er{prD} gives that the
function $\gf$ is real on the real line. \BBox

\section{\lb{Sec6} Appendix: Inverse problems for finite Jacobi operators}
\setcounter{equation}{0}

\subsection {Classis of Jacobi operators}

We define the class of Jacobi operators $H=H_o+V$ with triangular
matrix coefficients:
\[
\begin{aligned}
\lb{ta} &\M_d^\na=\{a\in \M_d^o \ \ {\rm is\ a\ upper\ triangular\
matrix\ with\ diagonal\ positive\ entries}\},
\\
& \J_q^\na=\{V\in \J_q^c:   s_x^*=a_x\in \M_d^\na,  \qq b_x=b_x^*,
\qq \forall \ x\ge 1\},\qq  \M_d^u=\{c\in \M_d^o\ \ \  {\rm is\
unitary}\}.
\end{aligned}
\]
We discus presentations of Jacobi matrices.

\begin{lemma}\lb{Tuu}
Consider a Jacobi operator $H_c=H_o+V_c$, where $s_x^*=a_x\in
\M_d^o$ and a sequence $(a_x, b_x)\in \M_d^o\ts \M_d, x\ge 1$
satisfies $\sup_{x\ge 1} (|a_x|+|b_x|)<\iy$, where $|\cdot|$ is the
norm in $\C^d$. Introduce the unitary operator $U$ on $\mH$ by
$(Uf)_x=u_xf_x, x\in \N$. Then the operator $T=U^* H U=H_o+V$, where
the operator $V$ and
\[
\lb{U2}
\begin{aligned}
& (Vf)_x=(\ga_{x-1}-\1)y_{x-1}+(\ga_x^*-\1)y_{x+1}+\gb_xy_x,
\\
&\ga_x= u_{x+1}^* a_{x}^+ u_{x},\qq \qqq \gb_x=u_x b_xu_x^*,\qq x\in
\N.
\end{aligned}
\]
\no  i) If the sequence $u_x$ has the form
\[
\lb{u1}
\begin{aligned}
u_{x+1}=\o_{x}u_1,\qq \o_{x}=v_{x-1}\cdots v_2 v_1,\qq x\in \N
\end{aligned}
\]
where $v_x\in \M^u$ is uniquely defined by the identity $a_x=v_x
a_x^+, a_x^+>0$, then
\[
\ga_x= u_{x+1}^* a_{x}^+ u_{x}>0.
\]
\no  ii) Let the sequence $u_x\in \M^u$ be defined by the rule: \ \

\no  $\bu $ $u_1\in \M^u$ is any matrix and

\no  $\bu $  if $u_x, x\ge 1$ is given, then $u_{x+1}$ is a unique
matrix under the condition $u_{x+1}^*a_xu_x\in \M^\na$.

Then
\[
\ga_x=u_{x+1}^*a_xu_x\in \M_d^\na
\]
iii) In the both cases, we have that   $V_c\in \J^c$ iff $V\in
\J^c$.

\end{lemma}
{\no\bf Proof.}  We consider the case $\gb=(\gb_x)=0$, the proof for
$\gb\ne 0$ is similar. Introduce the operator $U$ by
$y_x=(Uf)_x=u_xf_x$ where $f=(f_x)_1^\iy$. For finitely supported
sequences $y=(y_x)_1^\iy$ we have the scalar products
\[
\lb{hy}
\begin{aligned}
(Ty,y)= \sum_{1}^\iy \big((a_x y_x, y_{x+1})+(a_x^*y_{x+1},y_x)
\big)
  =\sum_{1}^\iy \big(
(\ga_xf_x,f_{x+1})+ (\ga_x^*f_{x+1}, f_{x}) \big),
\end{aligned}
\]
where $\ga_x=u_{x+1}^*a_xu_x$.

\no i) Due to the identity $a_x=a_x^+v_x$, where $v_x$ is uniquely
defined, we have $ \ga_x=u_{x+1}^* a_x^+v_xu_{x}$. Assume that $u_x,
x\ge 1$ satisfy
\[
\lb{euv} u_{x+1}=v_xu_x , \qq x\ge 1,
\]
then $ \ga_x=u_{x+1}^* a_x^+v_x u_{x}=u_{x+1}^* a_x^+ u_{x+1}>0. $
We show that the equation \er{euv} has the solution. Let
$\o_x=v_x\cdots v_2 v_1$ for any $x\ge 2$ and $u_1$. Then we have:
$$
\begin{aligned}
u_2=v_1u_1,\qq u_3=v_2u_2=v_2v_1u_1=\o_2u_1,\qq \dots,
u_{x+1}=\o_xu_1,....
\end{aligned}
$$
Thus this sequence $u_x$ satisfies the equation \er{euv} for any
$u_1$ and in particular for $u_1=\1$. The uniqueness of $u_x$
follows from the uniqueness in  $a_x=a_x^+v_x$.

\no ii) We need the well known fact:

\no $\bu $ {\it Let a matrix $A\in \M_d^o$. Then there exist unique
matrices $L\in \M_d^\na$ and   $\cU\in \M_d^u$  such that}
\[
\lb{au}   \qqq \cU\ A =L.
\]
We take any $u_{1}\in \M^u$. Due to \er{au}  we have
$\ga_1=u_{2}^*a_1u_1\in \M^\na$ for some unique $u_{2}\in \M^u$.
Again due to \er{au} (2) we have $\ga_2=u_{3}^*a_2u_2\in \M^\na$ for
some unique $u_{3}\in \M^u$ and so on.  The uniqueness of $u_x$
follows from the uniqueness in \er{au}.

\no iii) We show that  $V_c\in \J^c_q$ iff $V\in \J^c_q$ for the
case i). The proof for ii) is similar. Thus $a_{x}=\1, v_{x}$ for
all $x>p$ and $u_{p+1}=u_{p+2}=...$. This yields
$\ga_{p+j}=u_{p+j+1}^*a_{p+j}u_{p+j}=u_{p+1}^*a_{p+j}u_{p+j}=\1$ for
all $j\ge 1$. For $x=p$ we have
$$
\det [\1-\ga_{p}^2]=\det
[\1-u_{p+1}^*a_pu_p(u_{p}^*a_p^*u_{p+1})]=\det [\1-a_pa_p^*], \qq
\det \gb_p=\det (u_p b_pu_p^*)=\det b_p,
$$
which yields that $V_c\in \J^c_q$ iff $V\in \J^c_q$.
  \BBox

The unitary operator $U$ in \er{U2} has the matrix $u_1\in \M_d^u$
as a free parameter. All other matrices are defined uniquely. If
$u_1=\1$, then the unitary operator $U$ is uniquely defined.

\subsection {Definitions }

Recall results from \cite{K21a}. We discuss a finite Jacobi operator
$\cJ_{p}$ acting on $y=(y_x)_1^p\in \mH_p=(\C^d)^p$ with the
Dirichlet boundary condition and defined by
 \[
 \label{Jp}
 \begin{aligned}
\ca (\cJ_p y)_x= a_xy_{x+1}+b_xy_x+a_{x-1}y_{x-1},\qqq
 x\in \N_p=\{1,2,\dots,p\},\\
a_0=a_p=\1_d,\qqq y_0=y_{p+1}=0, \ac,
\end{aligned}
\]
where $b_x=b_x^*, a_x>0$ for all $x$. We can consider $\cJ_p$ as a
$p\ts p$ self-adjoint matrix given by
\[
\label{Ja} \cJ_p=\left(\begin{array}{ccccccc} b_1 & a_1 & 0 & 0 &
... &  0 \cr a_1 & b_2 & a_2 & 0 &  ... &  0 \cr 0 & a_2 & b_3 & a_3
&  ... & 0 \cr ... & ... & ... & ... & ... & ...   \cr  0 & ... & 0
& a_{p-2} & b_{p-1} & a_{p-1}\cr 0 & ...  & 0 & 0 & a_{p-1} & b_p
\end{array}\right),\qqq  b_x=b_x^*,\qqq a_x>0.
\]
The spectrum of $\cJ_p$ is a finite sequence of real eigenvalues
$\m_j$, $j\in \N_m$ with multiplicity $d_j\le d$, which satisfy
$$
\begin{aligned}
\m_1<\m_2<\dots<\m_m, \qqq d_j\le d,\qq \sum_1^m d_j=dp.
 \end{aligned}
$$
Define the $\M^d$-valued solutions $\vp(\l)=(\vp_x(\l))_{0}^{p+1}$
and $\c(\l)=(\c_x(\l))_{0}^{p+1}$ of the equation
\[
 \label{001}
  \begin{aligned}
& a_xy_{x+1}+b_xy_x+a_{x-1}y_{x-1}=\l y_x,\qqq
 \ \ x=0,1,..,p+1\,,
\\
& \vp_{0}=\c_{p+1}=0,\qq \vp_1=\c_{p}=\1_d,\qqq \l\in \C,
  \end{aligned}
\]
%for all $$.
Recall that $\vp_x(\l)$ and $\c_x(\l)$, $x\in\N_p$ are
matrix-valued polynomials such that $ \deg\vp_x(\l)=x$ and
$\deg\c_x(\l)=p\!+\!1\!-\!x. $
 If $\cJ_p y=\m_j y$, then the eigenvector $y\in \mH_p$
corresponding to the eigenvalue $\m_j$ has the form $y=(y_x)_1^p$,
where $y_x=\vp_x(\m_j)v$ for some $v\in\C^d$. Hence, the eigenvalues
of $\cJ_p$ are zeros of the polynomial $\det \vp_{p+1}(\cdot)$. We
sometimes write $\vp_x(\l,J_p), \m_{j}(J_p),...$ instead of
$\vp_x(\l), \m_{j}, ...$, when several Jacobi matrices are being
dealt with.

\no {\bf Definition of Spectral data.} {\it For each eigenvalue
$\m_j, j\in \N_m$, we define {the subspace}
\[
\label{cEd} {\cE_j}=\Ker \vp_{p+1}(\m_j),\qq \dim {\cE_j}=d_j\le d.
\]
$\bu$ {\bf the orthogonal projector} $\cP_j:\C^d\to\cE_j$ onto
$\cE_j\ss\C^d$,

\no $\bu$ {\bf the positive self-adjoint operator}
$g_j:\cE_j\to\cE_j$ given by}
\[
\label{Gdef} {g_j}= G_j|_{\cE_j},\qq {where}\qq
G_j=\cP_j\sum_{n=1}^{p} \vp_x^*(\m_j) \vp_x(\m_j)\cP_j.
\]
We  define  the matrix-valued Weyl-Titchmarsh function $M(\l)$ by
\[
\label{mdef} M(\l)=-\c_1(\l)\c_0(\l)^{-1},\qqq \l\in\C.
\]
 Define the reflection operator $\flat_p$ acting on
$y=(y_x)_1^p\in (\C^p)^m$ by
\[
\lb{dRef} (\flat_p y)_j=y_{p+1-j}, \qq \forall \ \ j\in \N_p,
\]
and the Jacobi matrix $\cJ_p^{\flat}:=\flat_p \cJ_p \flat_p$ under
the reflection. Consider the relationship between our spectral data
and Weyl-Titchmarsh function.

\begin{proposition}
\label{Twt} i) The fundamental solutions $\c_0, \c_1$ satisfy
\[
\begin{aligned}
 \lb{c1x}
 \c_0(\l,J_p)=\vp_{p+1}^*(\l,J_p)=\vp_{p+1}(\l,J_p^{\flat}),\\
 \c_1(\l,J_p)=-\vt_{p+1}^*(\l,J_p)=\vp_{p}(\l,J_p^{\flat}).
\end{aligned}
\]
\no ii) The function $M(\l)=M^*(\ol{\l})$ is analytic in $\C\sm
\s(\cJ_p)$ and has the form
\[
\label{iwt} M(\l)=-\sum_{j=1}^m {B_j\/\l-\m_j},\qqq\qqq \sum_{j=1}^m
B_j=\1_d,
\]
where the matrices $B_j=\res_{\l=\m_j}M(\l)$ are self-adjoint and
are given by
\[
\label{Bdef} B_j\big|_{\cE_j}=g_j^{-1},\qqq
B_j\big|_{\C^q\om\cE_j}=0,\qq j\in\N_m.
\]
iii) An operator $\mT$ on $(\C^d)^p$ given by \er{mT} is positive,
i.e.,
\[
\label{mT} \mT= \left(\begin{array}{cccc} \dP_0 & \dP_1 & ... &
\dP_{p-1} \cr \dP_1 & \dP_2 & ... & \dP_p \cr ... & ... & ... & ...
\cr \dP_{p-1} & \dP_{p} & ... & \dP_{2p-2}\end{array}\right)>0,\qq
{\rm where}\qq \dP_s=\sum_{j=1}^m\m_j^s\cP_j,\qq s=0,..,{2p\!-\!2}.
\]

\end{proposition}

In order to formulate our main result we need a restriction on
spectral data.

\noindent {\bf Example.} 1) Let $p=m$ and $d_j=\dim \cE_j = \rank
\cP_j = d$ for all $j\in\N_m$, i.e., $\cP_j=\1$ for all $j$. Then
for each distinct values $\m_1<\dots<\m_m$ the system
$(\m_j,\1)_1^m$ satisfies \er{mT}.

\no 2) Let $p=2$, $m=3$, $d_1+d_2=d_3=d$. Then the system
$\{(\m_1,\cP_1),(\m_2,\cP_2),(\m_3,\1)\}$ satisfies \er{mT} iff
$\Ker \cP_1\cap\Ker \cP_2=\{0\}$. Thus, in the continuous case
terminology, $\Ker \cP_2$ is the "forbidden subspace"\ (see
\cite{CK06}) for the projector $\cP_1$ and vice versa.

\medskip

For each $d\ge 2, p\ge 2$ we introduce the set of spectral data
\[
\label{Sdef} \gD_p= \left\{ \left(\m_j,\cP_j,g_j\right)_1^m:
\begin{array}{l} p\le m\le pd,\ \m_1<\dots<\m_m {\rm \ are\ real\
numbers}; \cr (\m_j,\cP_j)_1^m\ {\rm satisfies\ \er{mT}} \cr
g_j:\Ran \cP_j\to\Ran \cP_j\ {\rm are\ linear\ operators\ such} \cr
{\rm that\ } g_j=g_j^*>0\ {\rm and\ } \sum_{j=1}^m
\cP_jg_j^{-1}\cP_j =\1_d
\end{array}\right\}.
\]
We formulate our main result for the mapping $\P$ given by
$
\Psi: \left((a_x)_1^{p-1};(b_x)_1^p\right)\mapsto
\left(\m_j,\cP_j,g_j\right)_1^m\,.
$

\begin{theorem}\label{Thm}
 The mapping $\Psi:(\M_d^+)^{p-1}\!\ts (\M_d^s)^p\to\gD_p$ is
one-to-one and onto.

\end{theorem}

\no {\bf Acknowledgments.} \footnotesize   Our study was supported
by the RSF grant No 19-71-30002.

\end{document}